\documentclass{svjour3}
\usepackage{amssymb}

\usepackage{epic}
\usepackage{amsmath}
\usepackage{amsfonts}
\usepackage{pstricks}

\usepackage{float}
\usepackage{subfig}
\usepackage{graphicx}
\usepackage{multirow}
\usepackage{adjustbox}
\usepackage{color}
\usepackage[subfigure]{tocloft}
\usepackage{subfig}
\usepackage[utf8]{inputenc}
\usepackage{soul}

\allowdisplaybreaks
\numberwithin{equation}{section}
\usepackage[normalem]{ulem}
\usepackage[hidelinks]{hyperref}

\newenvironment{rcases*}{\left.\begin{array}{ll}}{\end{array}\right\rbrace}

\def\bigO{\mathcal{O}}

\def\L{\mathcal{L}}
\def\sii{\Leftrightarrow}
\newcommand{\C}{\mathcal{C}}

\def\bigO{\mathcal{O}}
\journalname{Journal of Scientific Computing}

\begin{document}
\title{On the efficient computation of smoothness indicators for a class of  WENO reconstructions}
\titlerunning{Efficient smoothness indicators for a class of WENO methods}
\author{Antonio~Baeza \and Raimund~B\"{u}rger \and Pep~Mulet \and David~Zor\'{\i}o}
\authorrunning{A.\ Baeza, R.\ B\"{u}rger,  P.\ Mulet and D.\ Zor\'{\i}o}
\date{}

\institute{A.\ Baeza \and P.\ Mulet  \at 
  Departament de Matem\`{a}tiques \\  Universitat de
  Val\`{e}ncia \\ E-46100 Burjassot,      Spain \\ 
  \email{antonio.baeza@uv.es, mulet@uv.es} 
  \and
  R.\ B\"{u}rger \at  
  CI$^2$MA \& Departamento de Ingenier\'{\i}a Matem\'{a}tica \\ 
   Universidad de Concepci\'{o}n  \\ 
   Casilla 160-C,  Concepci\'{o}n,    Chile\\ 
   \email{rburger@ing-mat.udec.cl} 
   \and 
   D.\ Zor\'{\i}o {\bf
(Corresponding author) } \at 
    CI$^2$MA,  
   Universidad de Concepci\'{o}n  \\ 
   Casilla 160-C,  Concepci\'{o}n,    Chile\\ 
   \email{dzorio@ci2ma.udec.cl}  
}

\thispagestyle{empty}

\noindent This version of the article has been accepted for publication, after a peer-review process, and is subject to Springer Nature’s AM terms of use, but is not the Version of Record and does not reflect post-acceptance improvements, or any corrections. The Version of Record is available online at:

\noindent \url{https://doi.org/10.1007/s10915-019-00974-7}

\newpage
\setcounter{page}{1}

\maketitle

\begin{abstract}
  Common smoothness indicators used in Weighted Essentially Non\--Os\-cil\-la\-to\-ry (WENO) reconstructions [Jiang, G.S., Shu, C.W.: Efficient implementation of {Weighted} {ENO} schemes, J.\ Comput.\ Phys. \textbf{126}, 202--228 (1996)] have quadratic cost with respect to the order. A set of novel smoothness indicators with linear cost of computation with respect to the order is presented. These smoothness indicators can be used in the context of schemes of the type 
   introduced by Yamaleev and Carpenter  [Yamaleev, N.K., Carpenter, M.H.: A systematic methodology to for constructing high-order energy stable WENO schemes. J.\ Comput.\ Phys. \textbf{228}(11), 4248--4272 (2009)]. The accuracy properties of the resulting non-linear weights are the same 
    as  those arising from using the traditional Jiang-Shu smoothness indicators in Yamaleev-Carpenter-type reconstructions. The increase of the efficiency and ease of implementation are shown.

\end{abstract}

\section{Introduction}

\subsection{Scope} \label{subsec:scope}
Weighted Essentially Non-Oscillatory (WENO) schemes are a very useful and powerful tool to 
 reconstruct functions from discrete data. They avoid an oscillatory behaviour in presence of a discontinuity and attain the optimal interpolation order whenever possible when the data is smooth. 
  The most common contexts in which these schemes are used are  finite-difference or finite-volume schemes discretizing hyperbolic conservation laws, whose solutions are typically non-smooth (weak solutions) and usually present strong shocks. 
 WENO reconstructions are building blocks of numerical schemes that properly handle discontinuous solutions,
 and represent a popular way to construct high-order methods.
However, it is well known that despite the recognized efficiency of these schemes, the most expensive part of the algorithm is the computation of local  smoothness indicators, whose computational cost is quadratic with respect to the order of the scheme. 

It is the purpose 
 of this work  to advance a novel design of  smoothness indicators so that they are cheaper to compute (namely, at linear cost in terms 
 of the order of accuracy) than the Jiang-Shu smoothness indicators \cite{JiangShu96} in the context of 
  Yamaleev-Carpenter reconstructions \cite{YamaleevCarpenter2009}, while   both  properties of optimal accuracy   in case of smoothness and  robust  capture of discontinuities are ensured. To do so, we  first advance theoretical considerations that are needed for their  foundation, 
   and that  allow us to  ultimately define the new smoothness indicators along   with the associated weight design. The resulting 
    WENO schemes are addressed as ``Fast WENO'' (FWENO) schemes.

 \subsection{Related work} \label{subsec:related}
 The classical WENO schemes were proposed by Jiang and Shu \cite{JiangShu96,shu98} as an improvement of the original proposal of Liu et al. \cite{LiuOsherChan94}. The idea is to build a weighted combination of interpolators, with the weights depending on smoothness indicators that tune the weight according to the data. The smoothness indicators defined by Jiang and Shu are designed in a way such that they take small values if the data used to construct the indicator is smooth and large values otherwise, with the additional property that the values of indicators constructed with smooth data are close to each other. This property is important as the order of accuracy of the final reconstruction is strongly dependent on it.
 
  Several years later, Yamaleev and Carpenter \cite{YamaleevCarpenter2009O3} proposed a new third-order WENO scheme, which was later on extended to arbitrary order in \cite{YamaleevCarpenter2009} (henceforth, YC-WENO scheme). In this case, the non-linear weights are based instead on a ratio between a high-order undivided difference, which is very small when the data is smooth and large if a discontinuity crosses the data stencil, and the original Jiang-Shu smoothness indicators. In this latter case the accuracy of the scheme is based on the asymptotic convergence of the undivided difference and the smoothness indicators, rather than the closeness between the latter ones.  This allows one to simplify the smoothness indicators, which is presented in this paper.

\subsection{Outline of this paper}

The remainder of the paper is organized as follows.  In Section \ref{prel}, we provide the theoretical background for the novel smoothness indicators. 
 Section~\ref{sec:mod_weno} is devoted to the definition of the smoothness indicators, the construction of the corresponding numerical scheme and the analysis of its accuracy. Section~\ref{sec:num_exp}  contains several numerical experiments in which the new algorithm is compared against the previous ones  in terms of accuracy and efficiency. Finally, in Section \ref{sec:conclusions} some conclusions are drawn.

\section{Preliminaries}\label{prel}

\subsection{Asymptotic properties of functions} 
We recall that for $\alpha \in \mathbb{Z}$, 
\begin{align*}
  f(h)&=\bigO(h^\alpha)\sii \limsup_{h\to 0} \left|\frac{f(h)}{h^\alpha}\right| <\infty, 
  \end{align*} 
and define the more restrictive property 
\begin{align*}   
  f(h)&=\bar{\bigO}(h^\alpha) :\sii \limsup_{h\to 0} \left|\frac{f(h)}{h^\alpha}\right|  <\infty \land \liminf_{h\to 0} \left|\frac{f(h)}{h^\alpha}\right| > 0. 
\end{align*}
For positive functions~$f$ and~$g$, the properties
\begin{align*}
  \limsup_{h\to 0} f(h)g(h)&\leq \limsup_{h\to 0} f(h) \limsup_{h\to
    0}  g(h), \\
 \liminf_{h\to 0} f(h)g(h)&\geq \liminf_{h\to 0} f(h) \liminf_{h\to 0} g(h)  
\end{align*}
imply that for $\alpha, \beta \in \mathbb{Z}$, 
$\bigO(h^{\alpha}) \bigO(h^{\beta}) = \bigO(h^{\alpha+\beta})$ and 
$\bar\bigO(h^{\alpha}) \bar\bigO(h^{\beta})=\bar \bigO(h^{\alpha+\beta})$.
Here and in what follows it is always understood 
  that expressions of the form $\bigO(h^{\alpha})$,  $\bar\bigO(h^{\alpha})$ correspond to $h \to 0$. 
Similarly, taking into account that for $A_i\subseteq (0, \infty)$, $\sup_{i} \inf A_i = \inf_{i}\sup
A_i^{-1} $ where we define $A_i^{-1}:=\{ 1/x \colon x\in A_i\}$, 
it follows that for  a positive function~$f$,
\begin{align*}
  \liminf_{h\to 0} f(h) = \limsup_{h\to 0} f(h)^{-1},
\end{align*}
therefore, if $f$~is  positive, then 
$f(h)=\bar{\bigO}(h^\alpha)$ implies 
$f(h)^{-1}=\bar{\bigO}(h^{-\alpha})$.
 
\subsection{Point values and cell averages of smooth functions}  
Finite-difference and finite-volume schemes 
	for hyperbolic conservation laws are based on discretizations of the solution by means of point values and cell averages, respectively.
	 The following lemmas state
	the same result for both cases so that we can analyze them in a unified way. 

\begin{lemma}\label{lemprel}
  Assume that a function  $\varphi\in\C^{n+2}$ satisfies  $\smash{\varphi^{(k)}(0)=0}$ for 
$k=1,\dots,n$ and~$\smash{\varphi^{(n+1)}(0)\neq 0}$. Then
$\varphi(h)=\bar\bigO(h^{n+1})$.
\end{lemma}

\begin{proof}
  The $(n+1)$-th order Taylor expansion of $\varphi$ yields
  \begin{align*}
    \varphi(h) &= \frac{\varphi^{(n+1)}(0)}{(n+1)!}h^{n+1} +
    \bigO(h^{n+2}),
    \end{align*} 
     which implies that 
     \begin{align*} 
    \lim_{h\to 0}\frac{\varphi(h)}{h^{n+1}} &=
    \frac{\varphi^{(n+1)}(0)}{(n+1)!}\neq 0,
  \end{align*}
  which in turn means that $\varphi(h)=\bar\bigO(h^{n+1})$. \hfill  \rule{1.5ex}{1.5ex}
  \end{proof} 

\begin{lemma}\label{lemacc}
  Let $c, d, z\in\mathbb{R}$. Assume that
  \begin{align*}
    \begin{cases}
  \text{\em $f^{(k)}(z)=0$ for $k=1,\dots,n$, $f^{(n+1)}(z)\neq 0$, and $f\in\C^{n+2}$} &  \text{\em 
    if $c+d\neq 0$}, \\
  \text{\em $f^{(2k-1)}(z)=0$ for $k=1,\dots,m$, $f^{(2m+1)}(z)\neq 0$, and $f\in\C^{2m+2}$} &  \text{\em 
    if $c+d= 0$.}
\end{cases}  
\end{align*}
Then 
\begin{align} \label{eq:pep4}
  f(z+dh)-f(z+ch)&=\bar\bigO(h^\nu), \\
  \label{eq:pep5}
  \frac{1}{h}\left(
    \int_{z+(d-1/2)h}^{z+(d+1/2)h} f(x) \, \mathrm{d}x
    -\int_{z+(c-1/2)h}^{z+(c+1/2)h} f(x) \, \mathrm{d} x
  \right)
  &=\bar\bigO(h^\nu),
\end{align}
where 
\begin{align*}
\nu&=
  \begin{cases}
    n+1& \text{\em if $c+d\neq 0$,}  \\
    2m+1& \text{\em if $c+d=0$.}  
  \end{cases}
  \end{align*} 
\end{lemma}

\begin{proof}
Our purpose is to apply Lemma \ref{lemprel} to the difference 
$$\varphi(h)=f(z+dh)-f(z+ch)$$ to obtain \eqref{eq:pep4} and alternatively, to the expression 
 \begin{align} \label{eq:phih} 
  \varphi(h)=
  \frac{1}{h}\left(
    \int_{z+(d-1/2)h}^{z+(d+1/2)h} f(x) \, \mathrm{d} x
    -\int_{z+(c-1/2)h}^{z+(c+1/2)h} f(x) \, \mathrm{d} x
  \right)
\end{align}
to obtain \eqref{eq:pep5}.  
For \eqref{eq:pep4}, 
it follows directly that 
\begin{align*}
  \varphi^{(k)}(0)=\mu(k, c, d)f^{(k)}(z), \quad \text{where $\mu(k, c, d)=d^k-c^k$.} 
\end{align*}
Since $c\neq d$  there  clearly holds that  $\mu(k, c, d)=0$ if and only
if $k$ is even and $c+d=0$.

For \eqref{eq:pep5}, we obtain from \eqref{eq:phih} multiplied by~$h$  
\begin{align*}
    \bigl(h\varphi(h)\bigr)'(h)=&-\Bigl(d-\frac{1}{2}\bigr)f \left(z+\Bigl(d-\frac12 \Bigr)h  \right)+\Bigl(d+\frac{1}{2}\Bigr)f \left(z + \Bigl(d+\frac12 \Bigr)h \right) \\
&+\Bigl(c-\frac{1}{2}\Bigr)f \left( z+ \Bigl(c-\frac12 \Bigr)h \right)-\Bigl(c+\frac{1}{2}\Bigr)f \left( z+ \Bigl(c+\frac12 \Bigr)h \right)
\end{align*} 
   and for $k\geq 1$, 
   \begin{align*} 
    &\bigl(h\varphi(h) \bigr)^{(k+1)}(0)\\=&
    \biggl(-\Bigl(d-\frac12\Bigr)^{k+1}+ \Bigl(d+\frac12 \Bigr)^{k+1}
    + \Bigl(c-\frac12 \Bigr)^{k+1}- \Bigl(c+\frac12 \Bigr)^{k+1} \biggr) f^{(k)}(z). 
  \end{align*}
  On the other hand, the Leibniz formula for higher derivatives yields
  \begin{align*}
     \bigl(h\varphi(h) \bigr)^{(k+1)}(0)&=(k+1)\varphi^{(k)}(0). 
  \end{align*}
  Therefore, for $k\geq 1$ we get $\smash{\varphi^{(k)}(0)  =\mu(k, c, d) f^{(k)}(z)}$, where 
  \begin{align*} 
      \mu(k, c, d) & =\frac{1}{k+1}\biggl(- \Bigl(d-\frac12 \Bigr)^{k+1}+ \Bigl(d+\frac12 \Bigr)^{k+1}
    + \Bigl(c-\frac12 \Bigr)^{k+1}- \Bigl(c+\frac12 \Bigr)^{k+1}\biggr).
 \end{align*} 
  Since
  \begin{align*}
    \Bigl(d+\frac12 \Bigr)^{k+1}- \Bigl(d-\frac12 \Bigr)^{k+1} & =
    \sum_{l=0}^{k+1}\binom{k+1}{l}d^{k+1-l}\left(\frac{1}{2^l}-\frac{(-1)^l}{2^l}\right) \\ & =
    \sum_{s=0}^{\lfloor k/2\rfloor}\binom{k+1}{2s+1}d^{k-2s}\frac{1}{2^{2s}}, 
  \end{align*}    
  we obtain 
  \begin{align}\label{eq:pep2}
    \mu(k, c, d)=
    \frac{1}{k+1}
    \sum_{s=0}^{\lfloor k/2\rfloor}\binom{k+1}{2s+1}\frac{1}{2^{2s}}(d^{k-2s}-c^{k-2s}). 
  \end{align}
  Now, if $c=-d$, then $\mu(k, c, d)=0$ if $k$ is even and
    \begin{align*}
    \mu(k, -d, d)=
    \frac{1}{k+1}
    \sum_{s=0}^{\lfloor
      k/2\rfloor}\binom{k+1}{2s+1}\frac{1}{2^{2s-1}}d^{k-2s}\neq 0
  \end{align*}
   if $k$ is
  odd, 
  since the sign of all summands is the sign of~$d$.
  On the other hand, if $c+d\neq 0$, since $c\neq d$, then $|c| \neq
  |d|$, so all summands in \eqref{eq:pep2} have the
  same sign, which yields $\mu(k, c, d)\neq 0$ for any $k\geq 1$.
  Therefore $\mu(k, c, d) =0   $ if and only if $k$~is odd and $c+d=0$. 
  In both cases it follows from the definition of~$n$ and~$m$ in the
  assumptions that  
  \begin{align*}
    \begin{cases}
      \text{$\varphi^{(k)}(0)=0$ for $k=1,\dots,n$, $\varphi^{(n+1)}(0)\neq 0$,
      and $\varphi\in\C^{n+2}$}&  \text{if $c+d\neq 0$}, \\
      \text{$\varphi^{(k)}(0)=0$ for $k=1,\dots,2m$, $\varphi^{(2m+1)}(0)\neq 0$,
     and $\varphi\in\C^{2m+2}$}&  \text{if $c+d= 0$,}
\end{cases}  
\end{align*}
so Lemma~\ref{lemprel} yields the final result. \hfill  \rule{1.5ex}{1.5ex}
\end{proof}

\section{Modified smoothness indicators and new weight design}\label{sec:mod_weno}

In this section the modified WENO schemes with the new smoothness indicators are considered for YC-WENO-type schemes \cite{YamaleevCarpenter2009}.
Schemes of order $2r-1$ are based on   a stencil 
\begin{align} \label{Sdef} 
S=\{f_{-r+1},\ldots,f_{r-1}\},
\end{align}
 where $f_i$ results either from a point-value or a cell-average discretization of a function $f$: 
\begin{align*} 
f_i=\L[f](x_i) := \begin{cases} f(x_i) & \text{for point values,} \\
\displaystyle \frac{1}{h} \int_{x_i - h/2}^{x_i + h/2}  f( \xi) \, \mathrm{d} \xi & \text{for cell averages} \end{cases} 
\end{align*} 
for constant $h:=x_{i+1}-x_i$, where we wish to  reconstruct~$f$ at~$x_{1/2}$. In the following sections we introduce the 
``Fast WENO'' (FWENO) schemes.

\subsection{Fast WENO (FWENO) schemes of order $2r-1$, $r\geq2$}\label{sec:eweno}
The traditional Jiang-Shu smoothness indicators \cite{JiangShu96} are defined by
\begin{align}\label{eq:oldindicators} 
I_{r,i}^{\text{JS}}:=\sum_{k=1}^{r-1}\int_{x_{-1/2}}^{x_{1/2}}h^{2k-1}\big(p_{r,i}^{(k)}(x)\big)^2\mathrm{d}x,\quad 0\leq i\leq r-1,
\end{align} 
where $p_{r,i}$ are the corresponding interpolating polynomials associated to the substencils $S_{r,i}=\{f_{-r+1+i},\ldots,f_i\}$, $0\leq i\leq r-1$.
These smoothness indicators have been typically used in the literature
involving WENO schemes, although their evaluation is computationally
expensive (quadratic with respect to the order, as becomes evident from
the identity deduced from \cite[Proposition~5]{SINUM2011} (see
  subsection \ref{sec:efficiency} for a more efficient computation,
  but still quadratic in $r$):
\begin{align*} 
I_{r,i}^{\text{JS}}=\sum_{j=0}^{r-1}\sum_{k=0}^je_{j,k}f_{-r+1+i+j}f_{-r+1+i+k}.
\end{align*} 

We propose new smoothness indicators for both (point-value and cell-average) 
 reconstructions that have {\em linear}  cost with respect to the order (namely, 
 they involve  $r-2$ additions and $r-1$ multiplications), and that are defined by
\begin{align} \label{eq:newindicators} 
I_{r,i}:=\sum_{j=1}^{r-1}(f_{-r+i+j+1}-f_{-r+i+j})^2,\quad0\leq i\leq r-1.
\end{align} 
A detailed analysis on the computational cost of the smoothness indicators in both cases is included in Section \ref{sec:efficiency}.
The remaining parts of the algorithm are the same as  those defined in \cite{YamaleevCarpenter2009} for the  YC-WENO schemes. For the sake of exposition, we briefly describe it:

\smallskip 

\noindent Input: $S=\{f_{-r+1},\ldots,f_{r-1}\}$, with $f_i=\L[f](x_i)$, and $\varepsilon>0$.
\begin{enumerate}
\item Compute interpolating polynomials  $$p_{r,i}(x)=\mathcal{I}_m(x_{-r+1+i},\ldots,x_i;f_{-r+1+i},\dots,f_i; x), \quad 0\leq i\leq r-1,$$
where $\mathcal{I}_m$ computes approximate point values from either point values or cell averages, according to the discretization.

\item Compute the new smoothness indicators \eqref{eq:newindicators}. 

\item Obtain the  corresponding squared undivided differences of order $2r-2$:
  \begin{align} \label{difnodiv}
  d_r=\Biggl(\sum_{j=-r+1}^{r-1}(-1)^{j+r-1}\binom{2r-2}{j+r-1}f_j\Biggr)^2.
  \end{align}
\item Compute the terms
  \begin{align} \label{alphayc} 
  \alpha_{r,i}=c_{r,i}\biggl(1+\frac{d_r^{s_1}}{I_{r,i}^{s_1}+\varepsilon}\biggr)^{s_2},\quad0\leq i\leq r-1, 
  \end{align} 
  where  $c_{r,i}$ are the ideal linear weights, for some $s_1,s_2$ chosen by the user such that $s_1\geq 1$ and $s_2\geq\frac{r}{2s_1}$.
\item Generate the WENO weights:
  \begin{align} \label{omegari} 
  \omega_{r,i}=\frac{\alpha_{r,i}}{\alpha_{r,0} + \dots+ \alpha_{r,r-1}}, \quad 
   i =0, \dots, r-1. \end{align} 
\item Obtain the reconstruction at $x_{1/2}$:
  \begin{align} \label{qrx} 
  q_r(x_{1/2})=\sum_{i=0}^{r-1}\omega_{r,i}p_{r,i}(x_{1/2}). \end{align} 
\end{enumerate}
Output: $q_r(x_{1/2})$.

\begin{remark}
  As stated above, the difference between our proposed WENO method and YC-WENO is the usage of the new smoothness indicators $I_{r,i}$ as defined in \eqref{eq:newindicators} in the former case and the usage of the classical smoothness indicators $I_{r,i}^{\text{JS}}$ as defined in \eqref{eq:oldindicators} in the latter case.
  In turn, the difference between JS-WENO  and YC-WENO schemes is that in the former case the coefficients $\alpha_{r,i}$ are defined by 
  \begin{equation}\label{alphajs}
    \alpha_{r,i}=\frac{c_{r,i}}{(I_{r,i}^{\text{JS}}+\varepsilon)^s}, \quad s\geq r/2,
  \end{equation}
instead of  \eqref{alphayc}. 
\end{remark}

\subsection{Accuracy properties of FWENO schemes}

We next analyze the accuracy properties of the novel smoothness indicators.

\begin{lemma}\label{si}
  Let $r\geq3$, and a grid be defined by  $x_i=z+(c+i)h$ for  $-r+1\leq i\leq r-1$. Assume that 
    $S$ given by \eqref{Sdef} is a stencil such that $f_i=\L[f](x_i)$,
    with $f^{(k)}(z)=0$ for $1\leq k\leq n$ and $f^{(n+1)}(z)\neq0$,
    $f\in\C^{n+1}$, $n \in$ $\mathbb{N}\cup\{0\}$, and assume that there exists 
     $m \in$ $\mathbb{N}\cup\{0\}$ such that $f^{(2m+1)}(z)\neq0$. 
     Furthermore, assume that the quantities $I_{r,i}$ are given by \eqref{eq:newindicators}. 
     Then $\smash{I_{r,i}=\bar{\bigO}(h^{2n+2})}$.

  On the other hand, if $r\geq2$ and $f$ has a unique discontinuity located in  $(x_{-r+1},x_{r-1})$, then  there 
   exist indices $i_0,i_1$ with  $-r+1\leq i_0,i_1\leq r-1$ such that
    $\smash{I_{r,i_0} =\bar{\bigO}(1)}$ and 
    $\smash{I_{r,i_1}=\bigO(h^2)}$.
\end{lemma}

\begin{proof}
  If $c\not\in\{1/2+j\mid-r+1\leq j\leq r-2\}$, then by Lemma \ref{lemacc} there holds  
  \begin{align*} 
  f_{j+1}-f_j=\bar{\bigO}(h^{n+1})  \Leftrightarrow 
   (f_{j+1}-f_j)^2=\bar{\bigO}(h^{2n+2}) \quad \text{for all $j= -r+1, \dots,  r-2$.} 
   \end{align*} 
  Therefore, in particular,  one has
  	\begin{align*} 
  I_{r,i}  
   & =\sum_{j=1}^{r-1}\bar{\bigO}(h^{2n+2})=\bar{\bigO}(h^{2n+2})
   \quad \text{for all $i=-r+1, \dots, r-1$.} 
   \end{align*} 
  Let us now assume that $c= 1 /2+j_0$ for some $j_0$, $-r+1\leq j_0\leq r-2$. Then, again by Lemma \ref{lemacc}:
  \begin{align*}
    f_{j_0+1}-f_{j_0}&=\bar{\bigO}(h^{2m+1}); \\
    f_{j+1}-f_j&=\bar{\bigO}(h^{n+1}) \quad \text{for all $j= -r+1, \dots, r-2$, $j\neq j_0$.} 
  \end{align*}
  Thus, we have the following:
  If $j_0\not\in\{-r+i+1,\ldots,i-1\}$, then
  $$I_{r,i} 
   =\sum_{j=1}^{r-1}\bar{\bigO}(h^{2n+2})=\bar{\bigO}(h^{2n+2})$$
  Otherwise, that is if $j_0 \in\{-r+i+1,\ldots,i-1\}$, we get 
  \begin{align*}
    I_{r,i} 
    &=(f_{j_0+1}-f_{j_0})^2+\sum_{j=1,j\neq j_0+r-i}^{r-1}(f_{-r+i+j+1}-f_{-r+i+j})^2\\
    &=\bar{\bigO}(h^{4m+2})+\sum_{j=1,j\neq j_0+r-i}^{r-1}\bar{\bigO}(h^{2n+2})=\bar{\bigO}(h^{4m+2})+\bar{\bigO}(h^{2n+2})=\bar{\bigO}(h^{2n+2}), 
  \end{align*}
   since clearly by definition $2m+1\geq n+1$. 
   
  Finally, if $f$ has a unique discontinuity at $(x_{-r+1},x_{r-1})$, then there exists $j_1 \in \{-r+1, \dots,  r-2\}$ such that $f_{j_1+1}-f_{j_1}=\bar{\bigO}(1)$, whereas $f_{j+1}-f_j=\bigO(h)$ for $j_1 \neq j \in \{ -r+1, \dots,  r-2\}$. 
  Hence, if we select for instance $i_0=0$ and $i_1=r-1$ if $j_1\geq0$, or $i_0=r-1$ and $i_1=0$ if $j_1<0$, then clearly
  $\smash{I_{r,i_0}=\bar{\bigO}(1)}$ and $\smash{ I_{r,i_1}=\bigO(h^2)}$.  
  
  \hfill \rule{1.5ex}{1.5ex}
\end{proof}

\begin{remark}
  The case $r=2$ for the FWENO method is the same as in the original YC-WENO method, since in this case the proposed smoothness indicators are the same. In this case, the statement of Lemma~\ref{si}  does not hold in general, since one can have $I_{2,i}=\bigO(h^{4m+2})$ if $c=1/2+i-1$.
\end{remark}

\begin{theorem}
  Under the same conditions and notation as in Lemma~\ref{si}, with $r\geq3$ and dropping the role of $\varepsilon$, there holds
  \[q_r(x_{1/2})=
  \begin{cases}
    f(x_{1/2})+\bigO(h^{2r-1}) & \text{if $n<2r-3$,}  \\
    f(x_{1/2})+\bigO(h^{n+1}) & \text{if $n\geq2r-3$,}  \\
    f(x_{1/2})+\bigO(h^r) & \textnormal{if a discontinuity crosses the data.}
  \end{cases}
  \]
\end{theorem}
\begin{proof}
  We first assume smoothness with a critical point of order $n$. Then, by Lemma~\ref{si}, the new smoothness indicators satisfy
  $I_{r,i}=\bar\bigO(h^{2n+2})$. 
  We consider in first place the case $n<2r-3$. Then there holds for $\alpha_{r,i}$:
  \begin{align*}
    \alpha_{r,i}&=c_{r,i}\biggl(1+\frac{d_r^{s_1}}{I_{r,i}^{s_1}}\biggr)^{s_2}=c_{r,i}\biggl(1+\frac{\bigO(h^{(4r-4)s_1})}{\bar\bigO(h^{(2n+2)s_1})}\biggr)^{s_2}\\
    &=c_{r,i}\bigl(1+\bigO (h^{2(2r-3-n)s_1} ) \bigr)^{s_2}=c_{r,i}+\bigO (h^{2(2r-3-n)s_1}).
  \end{align*}
  Hence, the non-linear weights satisfy
  \begin{align*}
    \omega_{r,i}&=\frac{\alpha_{r,i}}{\sum_{j=0}^{r-1}\alpha_{r,j}}=\frac{c_{r,i}+\bigO(h^{2(2r-3-n)s_1})}{\sum_{j=0}^{r-1}\bigl(c_{r,i}+\bigO(h^{2(2r-3-n)s_1})\bigr)}=\frac{c_{r,i}+\bigO(h^{2(2r-3-n)s_1})}{\sum_{j=0}^{r-1}c_{r,i}+\bigO(h^{2(2r-3-n)s_1})}\\
    &=\frac{c_{r,i}+\bigO(h^{2(2r-3-n)s_1})}{1+\bigO(h^{2(2r-3-n)s_1})}=c_{r,i}+\bigO(h^{2(2r-3-n)s_1}).
  \end{align*}
  On the other hand,
  $\smash{p_{r,i}(x_{1/2})=f(x_{1/2})+\bigO(h^M)}$ 
  with $M:=\max\{r,n+1\}$. 
  Therefore, denoting by 
  \begin{align*} 
  p_r(x_{1/2}):=\sum_{i=0}^{r-1}c_{r,i}p_{r,i}(x_{1/2})
  \end{align*} 
   the  value at $x_{1/2}$ of the optimal $(2r-1)$-th order polynomial, we have
  \begin{align*}
    q_r(x_{1/2})-f(x_{1/2})&=q_r(x_{1/2})-p_r(x_{1/2})+p_r(x_{1/2})-f(x_{1/2})\\
    &=\sum_{i=0}^{r-1}\omega_{r,i}p_{r,i}(x_{1/2})-\sum_{i=0}^{r-1}c_{r,i}p_{r,i}(x_{1/2})+\bigl(p_r(x_{1/2})-f(x_{1/2})\bigr)\\
    &=\sum_{i=0}^{r-1}(\omega_{r,i}-c_{r,i})p_{r,i}(x_{1/2})+\bigO(h^{2r-1})\\
    &=\sum_{i=0}^{r-1}(\omega_{r,i}-c_{r,i})\bigl(p_{r,i}(x_{1/2})-f(x_{1/2})\bigr)+\bigO(h^{2r-1})\\
    &=\sum_{i=0}^{r-1}\bigO(h^{2(2r-3-n)s_1})\bigO(h^M)+\bigO(h^{2r-1})\\
    &=\bigO(h^{2(2r-3-n)s_1+M})+\bigO(h^{2r-1}).
  \end{align*}
  Taking into account that $s_1\geq1$, we have
  $$q_r(x_{1/2})-f(x_{1/2})=\bigO(h^{2(2r-3-n)+M})+\bigO(h^{2r-1}).$$
  We next analyze the exponent of the left summand, splitting the discussion into two cases.
  On one hand, if $n\leq r-1$, then $M=r$ and there holds
  $$2(2r-3-n)+M=2(2r-3-n)+r=5r-6-2n\geq 5r-6-2(r-1)=3r-4\geq 2r-1,$$
  where the last inequality holds since by assumption $r\geq3$.

  On the other hand, if $n\geq r-1$, then $M=n+1$ and then, since by assumption $n\leq 2r-4$,
  $$2(2r-3-n)+M=2(2r-3-n)+(n+1)=4r-5-n\geq 4r-5-(2r-4)=2r-1.$$
   Thus,   for the case $n<2r-3$ there holds
  $q_r(x_{1/2})-f(x_{1/2})=\bigO(h^{2r-1})$. 
  
  In second place, we assume now $n\geq 2r-3$, then, using that $\sum_{i=0}^{r-1}\omega_{r,i}=1$, there holds
  \begin{align*}
    q_r(x_{1/2})&=\sum_{i=0}^{r-1}\omega_{r,i}p_{r,i}(x_{1/2})=\sum_{i=0}^{r-1}\omega_{r,i}\bigl(f(x_{1/2})+\bigO(h^{n+1})\bigr)\\
    &=\sum_{i=0}^{r-1}\omega_{r,i}f(x_{1/2})+\bigO(h^{n+1})=f(x_{1/2})+\bigO(h^{n+1}).
  \end{align*}
  Finally, if a discontinuity crosses the data, assume that  $J_0$~is the set of indices associated to the substencils $S_{r,i}$, $i\in J_0$, in which the discontinuity is not crossed ($J_0\neq\varnothing$ by the second part of Lemma~\ref{si}), and $J_1$ the set of the remaining whose corresponding substencils $S_{r,i}$, $i\in J_1$, are crossed by the discontinuity. Then if $i\in J_0$, 
  \begin{align*}
    \alpha_i=c_{r,i}\biggl(1+\frac{d_r^{s_1}}{I_{r,i}^{s_1}}\biggr)^{s_2}=c_{r,i}\left(1+\frac{\bar\bigO(1)}{\bar\bigO(h^{2s_1(n+1)})}\right)^{s_2}=\bar\bigO(h^{-2s_1s_2(n+1)}), 
  \end{align*}
  and  if $i\in J_1$, 
  \begin{align*}
    \alpha_i=c_{r,i}\biggl(1+\frac{d_r^{s_1}}{I_{r,i}^{s_1}}\biggr)^{s_2}=c_{r,i}\left(1+\frac{\bar\bigO(1)}{\bar\bigO(1)}\right)^{s_2}=\bar\bigO(1).
  \end{align*}
  Therefore, if $i\in J_0$, 
  \begin{align*}
    \omega_i&=\frac{\alpha_i}{\sum_{j=0}^{r-1}\alpha_j}=\frac{\alpha_i}{\sum_{j\in J_0}\alpha_j+\sum_{j\in J_1}\alpha_j}\\
    &=\frac{\bar\bigO(h^{-2s_1s_2(n+1)})}{\sum_{j\in J_0}\bar\bigO(h^{-2s_1s_2(n+1)})+\sum_{j\in J_1}\bar\bigO(1)}=\frac{\bar\bigO(h^{-2s_1s_2(n+1)})}{\bar\bigO(h^{-2s_1s_2(n+1)})}=\bigO(1),
  \end{align*}
  and if $i\in J_1$, 
  \begin{align*}
    \omega_i&=\frac{\alpha_i}{\sum_{j=0}^{r-1}\alpha_j}=\frac{\alpha_i}{\sum_{j\in J_0}\alpha_j+\sum_{j\in J_1}\alpha_j}\\
    &=\frac{\bar\bigO(1)}{\sum_{j\in J_0}\bar\bigO(h^{-2s_1s_2(n+1)})+\sum_{j\in J_1}\bar\bigO(1)}=\frac{\bar\bigO(1)}{\bar\bigO(h^{-2s_1s_2(n+1)})}=\bigO(h^{2s_1s_2(n+1)}).
  \end{align*}
  Thus, we conclude that 
  \begin{align*}
    q_r(x_{1/2})-f(x_{1/2})&=\sum_{i=0}^{r-1}\omega_{r,i}p_{r,i}(x_{1/2})-\sum_{i=0}^{r-1}\omega_{r,i}f(x_{1/2})\\
    &=\sum_{i=0}^{r-1}\omega_{r,i} \bigl( p_{r,i}(x_{1/2})-f(x_{1/2}) \bigr) \\
    &=\sum_{i\in J_0}\omega_{r,i} \bigl( p_{r,i}(x_{1/2})-f(x_{1/2})\bigr) +\sum_{i\in J_1}\omega_{r,i} 
     \bigl( p_{r,i}(x_{1/2})-f(x_{1/2}) \bigr)\\
    &=\sum_{i\in J_0}\bigO(1)\bigO(h^r)+\sum_{i\in J_1}\bigO(h^{2s_1s_2(n+1)})\bigO(1)\\
    &=\bigO(h^r)+\bigO(h^{2s_1s_2(n+1)}).
  \end{align*}
  Now using that $n\geq0$ and $s_2\geq\frac{r}{2s_1}$, we have
  $$2s_1s_2(n+1)\geq2s_1\frac{r}{2s_1}=r.$$
  Therefore $q_r(x_{1/2})-f(x_{1/2})=\bigO(h^r)$, 
  which completes the proof. \hfill  \rule{1.5ex}{1.5ex}
\end{proof}

\begin{remark}
  According to this result,   the only case in which our method loses accuracy with respect to the reconstruction with ideal linear weights is when $n=2r-3$, in which the accuracy order decays to $2r-2$, namely, one unit less than the optimal accuracy order, $2r-1$. 
\end{remark}

\subsection{Efficiency properties of FWENO schemes}
\label{sec:efficiency}
We conclude this section with a comparison involving the number of operations of an FWENO interpolator with respect to the traditional JS-WENO and YC-WENO interpolators. In order to do so, we invoke \cite[Propositions~1 and~5]{SINUM2011} to conclude that the evaluation of the reconstruction polynomials at the reconstruction point and the classical Jiang-Shu smoothness indicators can be respectively written as
\begin{align*} 
p_{r,i}(x_{1/2}) & =\sum_{j=0}^{r-1}d_{i,j}f_{-r+1+i+j}, \\
I_{r,i}^{\text{JS}} & =\sum_{j=0}^{r-1}\sum_{k=0}^je_{i,j,k}f_{-r+1+i+j}f_{-r+1+i+k},\quad 0\leq i\leq r-1,
\end{align*} 
where $d_{i,j},e_{i,j,k}\in\mathbb{R}$ are constants with respect to the data from the stencil.

This expression can be further simplified by taking into account that
$I_{r,i}^{\text{JS}}$ is a positive semi-definite quadratic form
defined on $(f_{-r+1+i},\dots,f_{i})$ with rank $r-1$, therefore it
can be expressed in a more convenient and numerically stable manner as
a sum of squares of linear combinations of
$f_{-r+1+i},\dots,f_{i}$. Specifically, for each $i=0,\dots,r-1$, let
$A_i$ be the matrix associated to $I_{r,i}^{\text{JS}}$, i.e.
\begin{align*}
  I_{r,i}^{\text{JS}} =
  (f_{-r+1+i},\dots,f_{i})
  A_i
  \begin{pmatrix}f_{-r+1+i}\\\vdots\\f_{i}  \end{pmatrix}.
\end{align*}
The $r\times r$ matrix
$A_i$ is  semi-positively-definite  with rank $r-1$
and therefore admits a decomposition as $A_i=P_i^{\mathrm{T}}L_i D_i L_i^{\mathrm{T}} P_i$,
where $P_i$ is a permutation matrix associated to the permutation
$\sigma_i$ of $(1,\dots,r)$, $L_i$ is lower triangular with
unit entries in the diagonal and
$D_i=\text{diag}(\beta_{i,1},\dots,\beta_{i,r-1}, 0)$, with $\beta_{i,j}>0$,
$j=1,\dots,r-1$.
This  yields the following simplified expression: 
\begin{align*}
  I_{r,i}^{\text{JS}}&=  (f_{-r+1+i},\dots,f_{i})
  P_i^{\mathrm{T}}L_i D_i L_i^{\mathrm{T}} P_i
  \begin{pmatrix}f_{-r+1+i}\\\vdots\\f_{i}  \end{pmatrix}\\
  &=
\sum_{j=1}^{r-1}\beta_{i,j}
\left(\sum_{k=0}^j\gamma_{i,j,k}f_{i-r+\sigma_i(1+k)}\right)^2,\quad
\gamma_{i,j,k}=(L_{i})_{j,k}, \gamma_{i,j,j}=1.
\end{align*}  

The cost associated to each reconstruction of lower order is thus
$r-1$ additions and $r$ multiplications, while the cost associated to
each classical smoothness indicators is
$(r^2+r-4)/2$ additions and
$(r^2+3r-4)/2$ multiplications. Hence, the cost
associated to the computation of the whole set of $r$ reconstructions
of lower order is $(r-1)r$ additions and $r^2$ multiplications,
whereas the cost associated to the computation of the whole set of $r$
classical smoothness indicators is
$r(r^2+r-4)/2$ additions and $r(r^2+3r-4)/2$ multiplications.

Now let us analyze the cost associated with the FWENO smoothness indicators \eqref{eq:newindicators}. In this case, we also can reduce the number of operations by taking into consideration that the novel smoothness indicators satisfy a simple recurrence relation, in a way that their computation can be simplified in the following algorithm, involving a linear cost both in terms of additions and in terms of multiplications:

\begin{enumerate}
\item Compute
  $$\theta_j:=(f_{-r+j+1}-f_{-r+j})^2,\quad1\leq j\leq 2r-2.$$
  Operations: $2r-2$ additions and $2r-2$ multiplications.
\item Compute
  $$I_0=\sum_{j=1}^{r-1}\theta_j.$$
  Operations: $r-2$ additions.
\item Compute
  $$I_i=I_{i-1}-\theta_i+\theta_{i+r-1},\quad 1\leq i\leq r-1.$$
  Operations: $2r-2$ additions.
\end{enumerate}

Therefore, the cost associated to the whole set of smoothness indicators involves $5r-6$ additions and $2r-2$ multiplications.

Now, YC-WENO and FWENO schemes also involve the computation of \eqref{difnodiv}, which involves $2r-2$ additions and $2r$ multiplications.

As for the terms $\alpha_{r,i}$, we have two cases: for JS-WENO schemes, the expression for them is \eqref{alphajs}, and therefore the number of operations for each one is one addition, $s-1$ multiplications and one division; therefore, the total cost for all them is $r$~additions, $(s-1)r$ multiplications and $r$~divisions. As for YC-WENO and FWENO schemes, the expression to be computed is \eqref{alphayc}, and the number of operations associated to each one of them is  two additions, $2s_1+s_2-2$ multiplications and one division, being thus the total cost $2r$ additions, $(2s_1+s_2-2)r$ multiplications and $r$ divisions.

The non-linear weights $\omega_{r,i}$ have the same expression  \eqref{omegari} in all cases. The denominator is the same for all the weights, and therefore one can previously store the value of $\smash{\bar{\alpha}_i:=1/(\sum_{j=0}^{r-1}\alpha_{i,j})}$ and then compute 
$\omega_{r,i}=\alpha_{r,i} \bar{\alpha}_i$, converting thus $r$~divisions, much more expensive than multiplications, 
 in to one division and $r$~multiplications. Therefore, the total cost corresponds now to $r-1$ additions and one  division associated 
  with  the computation of $\bar{\alpha}_{i}$ and $r$~multiplications (one multiplication per weight).
Finally, the reconstruction expression  \eqref{qrx} is also common in the three schemes, and corresponds to 
 $r-1$ additions and $r$~multiplications.  

The number of operations associated to each term and the grand total of operations for each method is summarized in Tables~\ref{effjsweno} to~\ref{effeweno}, where it can be drawn as a conclusion that the number of operations of JS-WENO and YC-WENO is cubic with respect to the order, whereas the number of operations associated to FWENO increases quadratically with respect to the order. Therefore, the complexity of the smoothness indicators is indeed a crucial factor in terms of the impact on the computational cost, and using simplified alternatives can indeed reduce significantly the overall computational cost of the WENO interpolator.

  
  

\begin{table}[t]
  \centering
  \begin{tabular}{|c|c|c|c|}
    \hline
    Order $2r-1$ & \multicolumn{3}{c|}{JS-WENO} \\
    \hline
    Operations & Additions & Multiplications & Divisions\\
    \hline
    $p_r(x_{\frac12})$ & $(r-1)r$ & $r^2$ & $0$\\
    \hline
    $I_r$ & $(r^3+r^2-4r)/2$ & $(r^3+3r^2-4r)/2$ & $0$\\
    \hline
    $\alpha_r$ & $r$ & $(s-1)r$ & $r$\\
    \hline
    $\omega_r$ & $r-1$ & $r$ & $1$\\
    \hline
    $q_r(x_{1/2})$ & $r-1$ & $r$ & $0$\\
    \hline
    SUM & $(r^3+3r^2-4)/2$ & $(r^3+5r^2+(2s-2)r)/2$ & $r+1$\\
    \hline
    TOTAL & \multicolumn{3}{c|}{$r^3+4r^2+(s-1)r-2$}\\
    \hline
  \end{tabular}
  
  \smallskip 
  
  \caption{Theoretical cost for JS-WENO interpolator.}
  \label{effjsweno}
\end{table}

\begin{table}[t]
  \centering
  \begin{tabular}{|c|c|c|c|}
    \hline
    Order $2r-1$ & \multicolumn{3}{c|}{YC-WENO} \\
    \hline
    Operations & Additions & Multiplications & Divisions\\
    \hline
    $p_r(x_{\frac12})$ & $(r-1)r$ & $r^2$ & $0$ \\
    \hline
    $I_r$ & $(r^3+r^2-4r)/2$ & $(r^3+3r^2-4r)/2$ & $0$ \\
    \hline
    $d_r$ & $2r-2$ & $2r$ & $0$ \\
    \hline
    $\alpha_r$ & $2r$ & $(2s_1+s_2-2)r$ & $r$ \\
    \hline
    $\omega_r$ & $r-1$ & $r$ & $1$ \\
    \hline
    $q_r(x_{1/2})$ & $r-1$ & $r$ & $0$ \\
    \hline
    SUM & $(r^3+3r^2+6r-8)/2$ & $(r^3+5r^2+(4s_1+2s_2)r)/2$ & $r+1$ \\
    \hline
    TOTAL & \multicolumn{3}{c|}{$r^3+4r^2+(2s_1+s_2+3)r-4$} \\
    \hline
  \end{tabular}
  
  \smallskip 
  
  \caption{Theoretical cost for YC-WENO interpolator.}
  \label{effycweno}
\end{table}

\begin{table}[t]
  \centering
  \begin{tabular}{|c|c|c|c|}
    \hline
    Order $2r-1$ & \multicolumn{3}{c|}{FWENO} \\
    \hline
    Operations & Additions & Multiplications & Divisions\\
    \hline
    $p_r(x_{\frac12})$ & $(r-1)r$ & $r^2$ & $0$ \\
    \hline
    $I_r$ & $5r-6$ & $2r-2$ & $0$ \\
    \hline
    $d_r$ & $2r-2$ & $2r$ & $0$ \\
    \hline
    $\alpha_r$ & $2r$ & $(2s_1+s_2-2)r$ & $r$ \\
    \hline
    $\omega_r$ & $r-1$ & $r$ & $1$ \\
    \hline
    $q_r(x_{1/2})$ & $r-1$ & $r$ & $0$ \\
    \hline
    SUM & $r^2+10r-10$ & $r^2+(2s_1+s_2+4)r-2$ & $r+1$ \\
    \hline
    TOTAL & \multicolumn{3}{c|}{$2r^2+(2s_1+s_2+14)r-12$} \\
    \hline
  \end{tabular}
  
  \smallskip 
  
  \caption{Theoretical cost for FWENO interpolator.}
  \label{effeweno}
\end{table}

\begin{figure}
  \centering
  \includegraphics[width=0.7\textwidth]{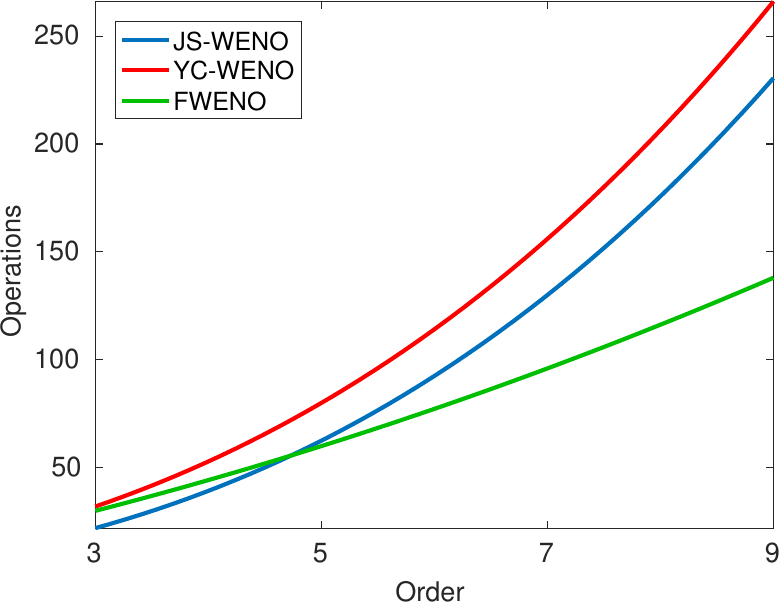}
  \caption{Theoretical cost comparison between WENO interpolators.}
  \label{effcomp}
\end{figure}

A graphical comparison between the cost associated to each scheme is also shown in Figure \ref{effcomp}, with $s=r/2$, $s_1=r/2$ and $s_2=1$.

\section{Numerical experiments}\label{sec:num_exp}

We now present some numerical experiments for schemes based on finite differences, as introduced in \cite{ShuOsher89,ShuOsher1989},  combined with high-order WENO reconstructions to discretize hyperbolic conservation laws.  Results obtained by the new FWENO scheme 
 are compared with those generated by the JS-WENO and YC-WENO schemes. The exponent $s$ for JS-WENO method is chosen as $s=\lceil r/2\rceil$, while the exponents $s_1$ and $s_2$ for both YC-WENO and FWENO methods are chosen as $s_1=\lceil r/2\rceil$ and $s_2=1$. Finally, we set in all the numerical experiments $\varepsilon=10^{-100}$, so that this parameter has the sole role of avoiding divisions by zero.

\subsection{1D conservation law experiments}\label{subsec:1d_cons_laws}

\subsubsection*{Example 1: Linear advection equation}

\begin{table}[t] 
  \centering {
  \begin{tabular}{|c|c|c|c|c|c|c|c|c|}
    \hline
    & \multicolumn{2}{c|}{$\|\cdot\|_1$} & \multicolumn{2}{c|}{$\|\cdot\|_{\infty}$} & \multicolumn{2}{c|}{$\|\cdot\|_1$} & \multicolumn{2}{c|}{$\|\cdot\|_{\infty}$} \\
    \hline
    & \multicolumn{4}{c|}{YC-WENO5} & \multicolumn{4}{c|}{FWENO5} \\
    \hline
    $N$ & Err. & $\bigO$ & Err. & $\bigO$ & Err. & $\bigO$ & Err. & $\bigO$ \\
    \hline
    10 & 1.02e-03 & --- & 1.55e-03 & --- & 1.01e-03 & --- & 1.66e-03 & --- \\
    20 & 3.27e-05 & 4.96 & 5.16e-05 & 4.91 & 3.27e-05 & 4.95 & 5.16e-05 & 5.01 \\
    40 & 1.01e-06 & 5.01 & 1.60e-06 & 5.01 & 1.01e-06 & 5.01 & 1.60e-06 & 5.01 \\
    80 & 3.15e-08 & 5.01 & 4.94e-08 & 5.01 & 3.15e-08 & 5.01 & 4.94e-08 & 5.01 \\
    160 & 9.79e-10 & 5.01 & 1.54e-09 & 5.01 & 9.79e-10 & 5.01 & 1.54e-09 & 5.01 \\
    320 & 3.05e-11 & 5.00 & 4.79e-11 & 5.00 & 3.05e-11 & 5.00 & 4.79e-11 & 5.00 \\
    640 & 9.52e-13 & 5.00 & 1.50e-12 & 5.00 & 9.53e-13 & 5.00 & 1.50e-12 & 5.00 \\
    \hline
  \end{tabular}}
  
  \smallskip 
  
  \caption{Example 1 (linear advection equation): fifth-order schemes with the new smoothness indicators.}
  \label{advlin_o5}
\end{table}

We first apply the fifth-order accurate  YC-WENO5 and FWENO schemes to 
 the initial-boundary value problem for the linear advection equation 
\begin{gather*} 
 u_t+ u_x=0,\quad\Omega=(-1,1),\quad u(-1, t)=u(1, t), \\
 u_0(x)=0.25+0.5\sin(\pi x),
 \end{gather*} 
 which 
 has the solution 
 $u(x,t)=0.25+0.5\sin(\pi(x-t))$. We run several simulations with final time $T=1$ and grid spacings $h=2/N$  and measure the resulting errors both 
  in  the $L^1$ and $L^{\infty}$   norms.  From the numerical results,  which  are shown in Table~\ref{advlin_o5}, 
 one can conclude that both schemes converge numerically at approximate fifth-order rate, which is consistent with our theoretical analysis. In fact, the results for  the YC-WENO and FWENO schemes are nearly identical.

\subsubsection*{Example 2: Burgers equation}

\begin{table}[t] 
  \centering {
  \begin{tabular}{|c|c|c|c|c|c|c|c|c|}
    \hline
    & \multicolumn{2}{c|}{$\|\cdot\|_1$} & \multicolumn{2}{c|}{$\|\cdot\|_{\infty}$} & \multicolumn{2}{c|}{$\|\cdot\|_1$} & \multicolumn{2}{c|}{$\|\cdot\|_{\infty}$} \\
    \hline
    & \multicolumn{4}{c|}{YC-WENO5} & \multicolumn{4}{c|}{FWENO5} \\
    \hline
    $N$ & Err. & $\bigO$ & Err. & $\bigO$ & Err. & $\bigO$ & Err. & $\bigO$ \\
    \hline
    40 & 2.44e-05 & 4.95 & 2.52e-04 & 4.70 & 2.73e-05 & 5.11 & 2.50e-04 & 4.69 \\
    80 & 7.89e-07 & 5.08 & 9.70e-06 & 5.05 & 7.89e-07 & 5.08 & 9.70e-06 & 5.05 \\
    160 & 2.32e-08 & 5.05 & 2.94e-07 & 5.03 & 2.32e-08 & 5.05 & 2.94e-07 & 5.03 \\
    320 & 7.01e-10 & 5.04 & 8.96e-09 & 5.03 & 7.01e-10 & 5.04 & 8.96e-09 & 5.03 \\
    640 & 2.14e-11 & 5.02 & 2.73e-10 & 5.02 & 2.14e-11 & 5.02 & 2.73e-10 & 5.02 \\
    1280 & 6.59e-13 & 5.02 & 8.41e-12 & 5.00 & 6.59e-13 & 5.02 & 8.41e-12 & 5.00 \\
    \hline
  \end{tabular}}
  
  \smallskip 
  
  \caption{Example 2 (Burgers equation): fifth-order schemes with the new smoothness indicators.}
  \label{burgers_o5}
\end{table} 

\begin{figure}[t] 
  \centering
  \begin{tabular}{cc}
    \multicolumn{2}{c}{\includegraphics[width=0.6\textwidth]{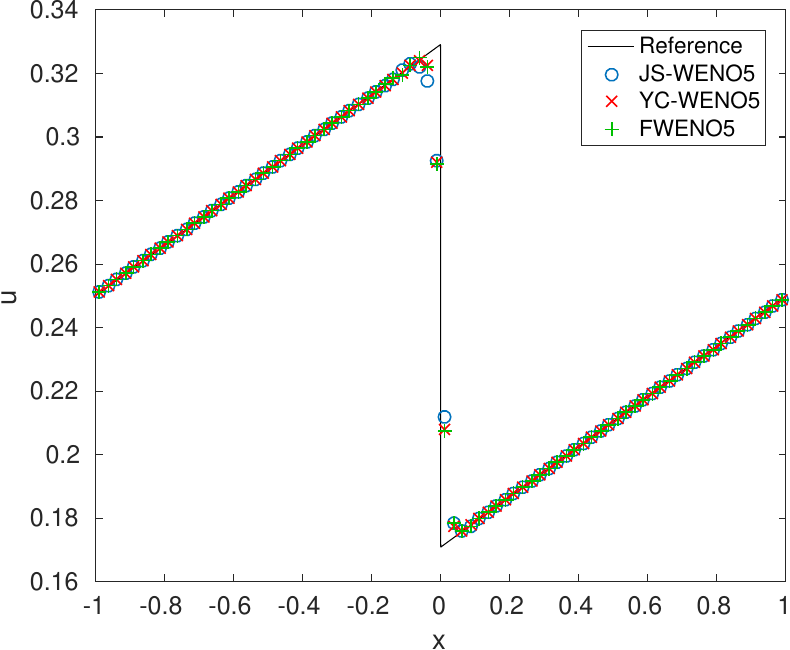}} \\
    \includegraphics[width=0.473\textwidth]{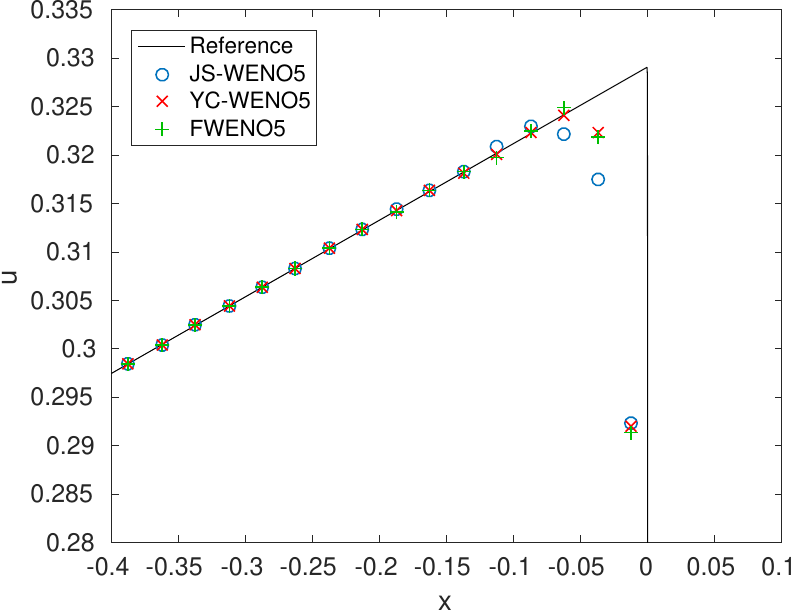} & \includegraphics[width=0.47\textwidth]{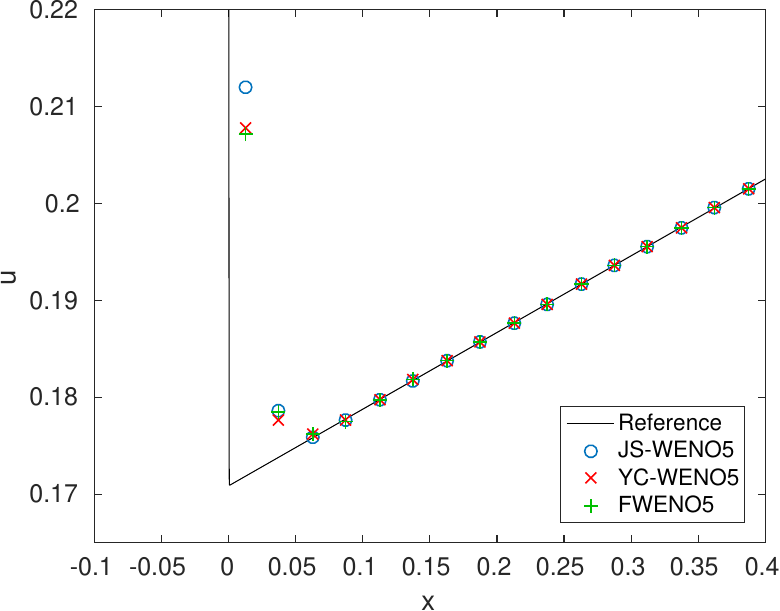} \\
  \end{tabular}
  \caption{Example 2 (inviscid Burgers equation): numerical solution for $N=80$ (top) and enlarged views (bottom left and right) at $T=12$ produced by 
   fifth-order schemes with conventional (JS-WENO5 and YC-WENO5) and 
   new (FWENO5) smoothness indicators. The reference solution with $N=16000$ cells is also shown.}
  \label{burgersdisc_o5}
\end{figure}

We now solve numerically the following initial-boundary value problem for the inviscid Burgers equation: 
\begin{gather*} 
u_t+(u^2/2)_x=0,\quad\Omega=(-1,1),\quad u(-1, t)=u(1, t), \\
 u_0(x)=0.25+0.5\sin(\pi x).
\end{gather*}

A first simulation is run until $T=0.3$, in which the solution remains smooth. We compare both YC-WENO and FWENO schemes in the results shown in Table \ref{burgers_o5}, where it can be observed that both methods attain again fifth order accuracy, producing an almost identical error.

The simulation is now run until $T=12$. At $t=1$, the wave breaks and a shock is generated. Therefore, we use the Donat-Marquina flux-splitting algorithm \cite{DonatMarquina96}. The results  shown in Figure~\ref{burgersdisc_o5} correspond to the fifth-order schemes, with a resolution of $N=80$ cells, and are compared with a reference solution computed with $N=16000$ cells. 
 The results obtained for all three methods are very similar, and therefore one can conclude that in this case 
  using the new smoothness indicators does not appreciably affect the  quality of the solution.

\subsubsection*{Example 3: Shu-Osher problem, 1D Euler equations of gas dynamics}

\begin{figure}[t] 
  \centering
  \begin{tabular}{cc}
    \includegraphics[width=0.47\textwidth]{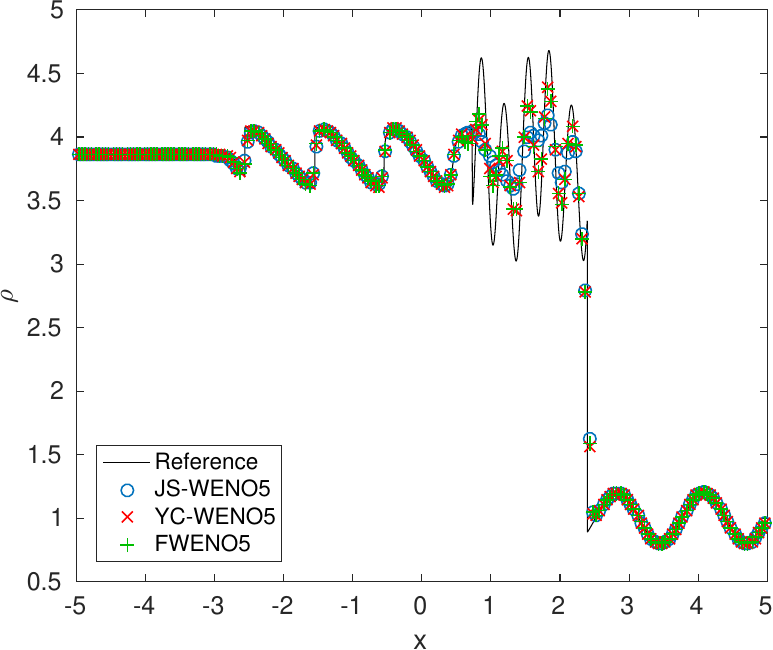} & \includegraphics[width=0.47\textwidth]{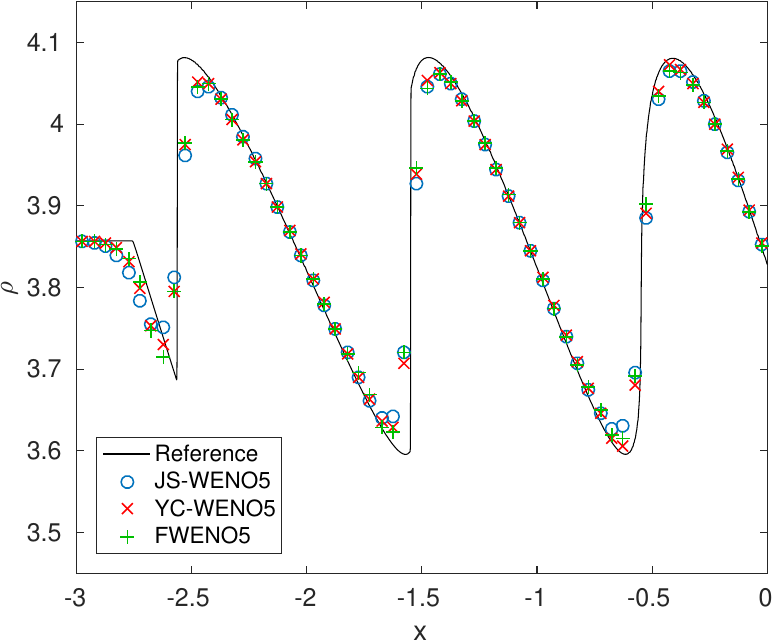} \\
    \includegraphics[width=0.47\textwidth]{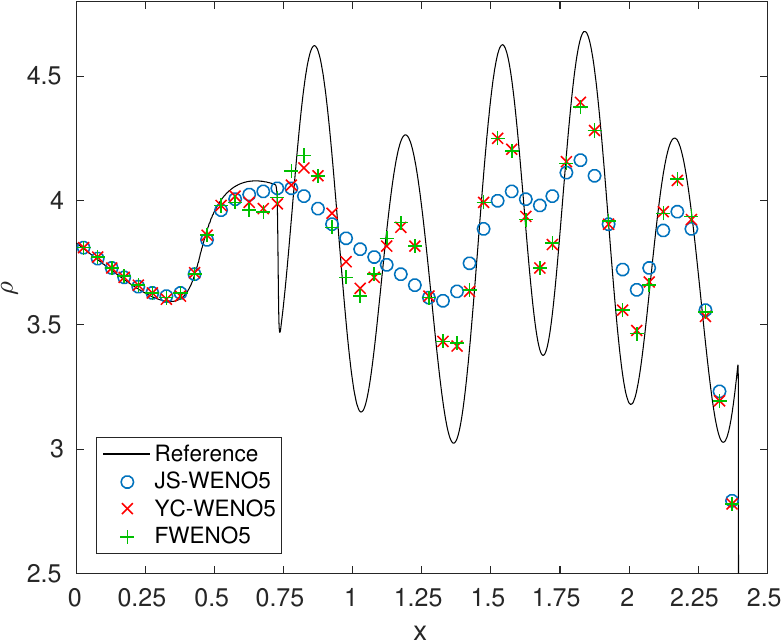} & \includegraphics[width=0.47\textwidth]{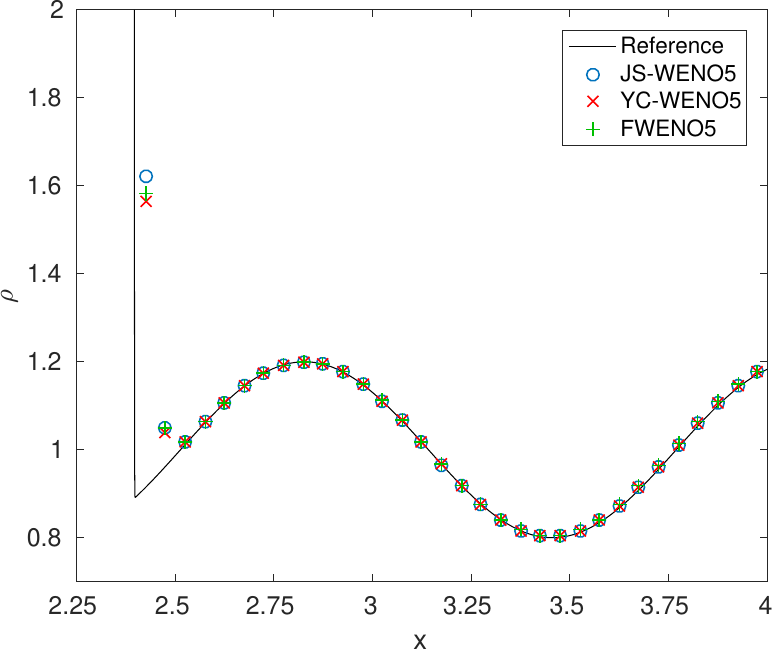}
  \end{tabular}
  \caption{Example 3 (Shu-Osher problem, 1D Euler equations of gas dynamics): numerical solution for $N=200$ (top left) and enlarged views (top right, bottom left and bottom right) at $T=1.8$ produced by 
   fifth-order schemes with conventional (JS-WENO5 and YC-WENO5) and 
   new (FWENO5) smoothness indicators. The reference solution with $N=16000$ cells is also shown.}
  \label{shuosher_o5_n=200}
\end{figure}

\begin{figure}[t]
  \centering
  \begin{tabular}{cc}
    \includegraphics[width=0.47\textwidth]{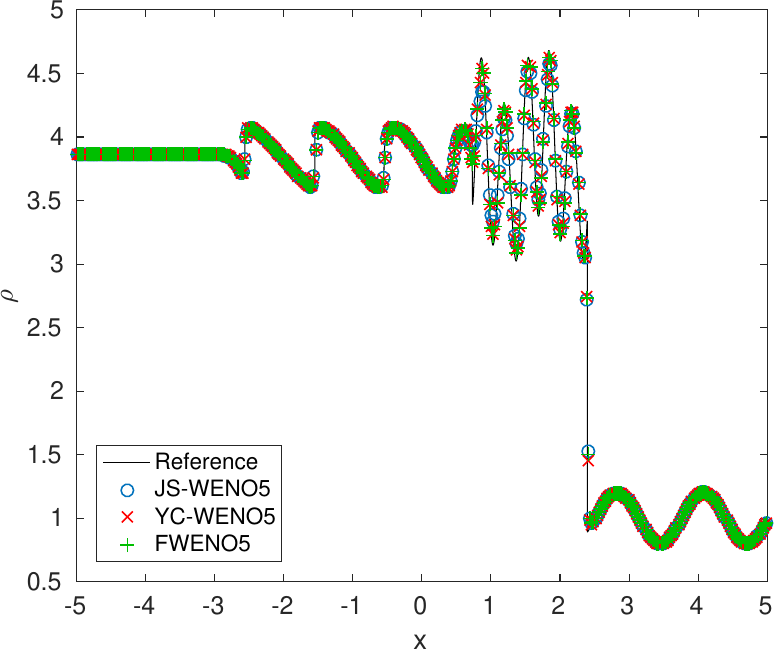} & \includegraphics[width=0.47\textwidth]{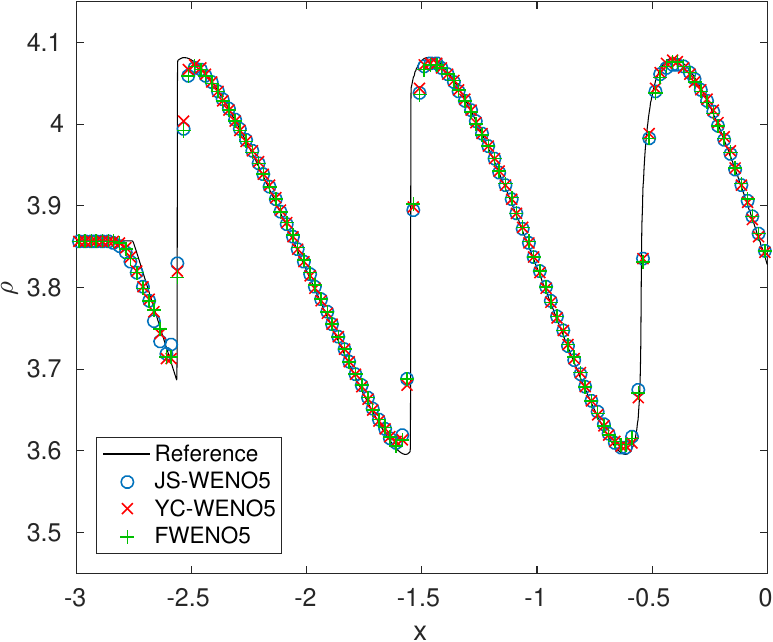} \\
    \includegraphics[width=0.47\textwidth]{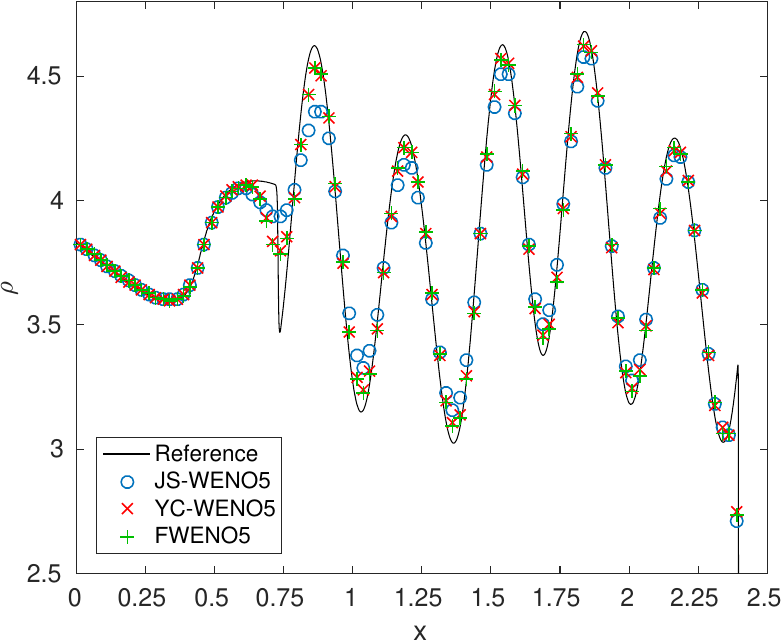} & \includegraphics[width=0.47\textwidth]{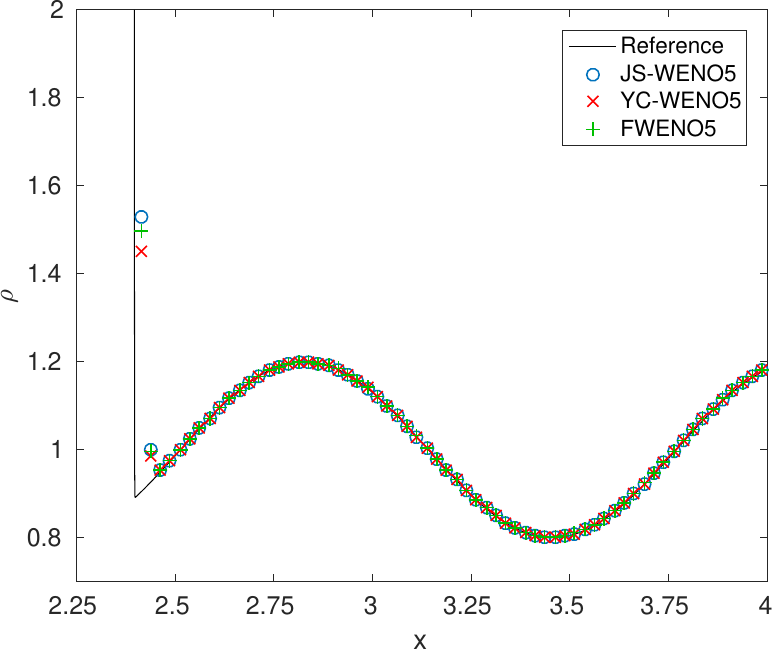}
  \end{tabular}
  \caption{Example 3 (Shu-Osher problem, 1D Euler equations of gas dynamics): numerical solution for $N=400$ (top left) and enlarged views (top right, bottom left and bottom right) at $T=1.8$ produced by 
   fifth-order schemes with conventional (JS-WENO5 and YC-WENO5) and 
   new (FWENO5) smoothness indicators. The reference solution with $N=16000$ cells is also shown. }
  \label{shuosher_o5_n=400}
\end{figure}

\begin{figure}[t]
  \centering
  \begin{tabular}{cc}
    \includegraphics[width=0.85\textwidth]{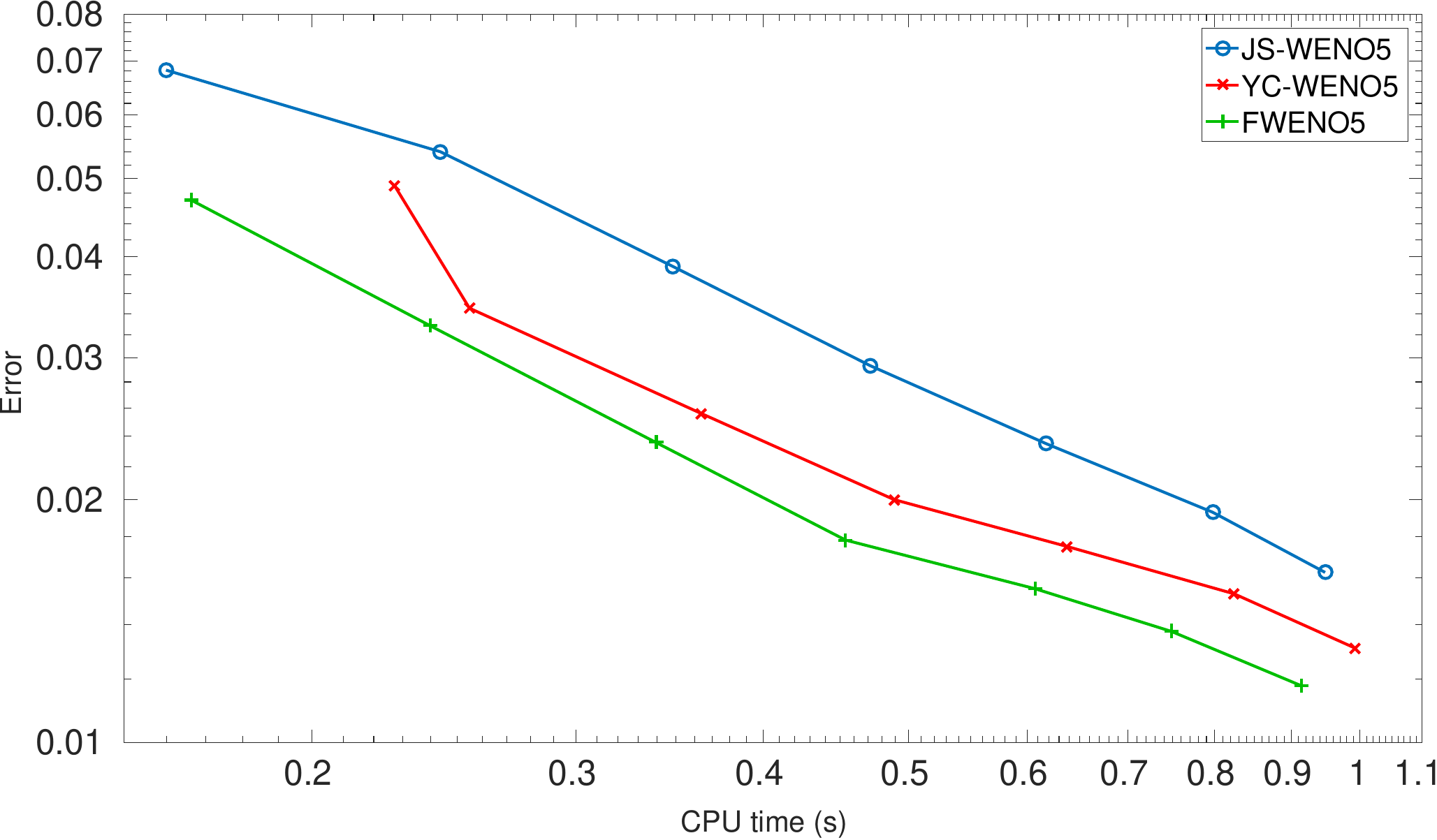}
  \end{tabular}
  \caption{Example 3 (Shu-Osher problem, 1D Euler equations of gas dynamics): efficiency plot for fifth-order schemes, corresponding to the 
   numerical solution at $T=1.8$.}
  \label{cpu_o5}
\end{figure}

\begin{figure}[t]
  \centering
  \begin{tabular}{cc}
    \includegraphics[width=0.85\textwidth]{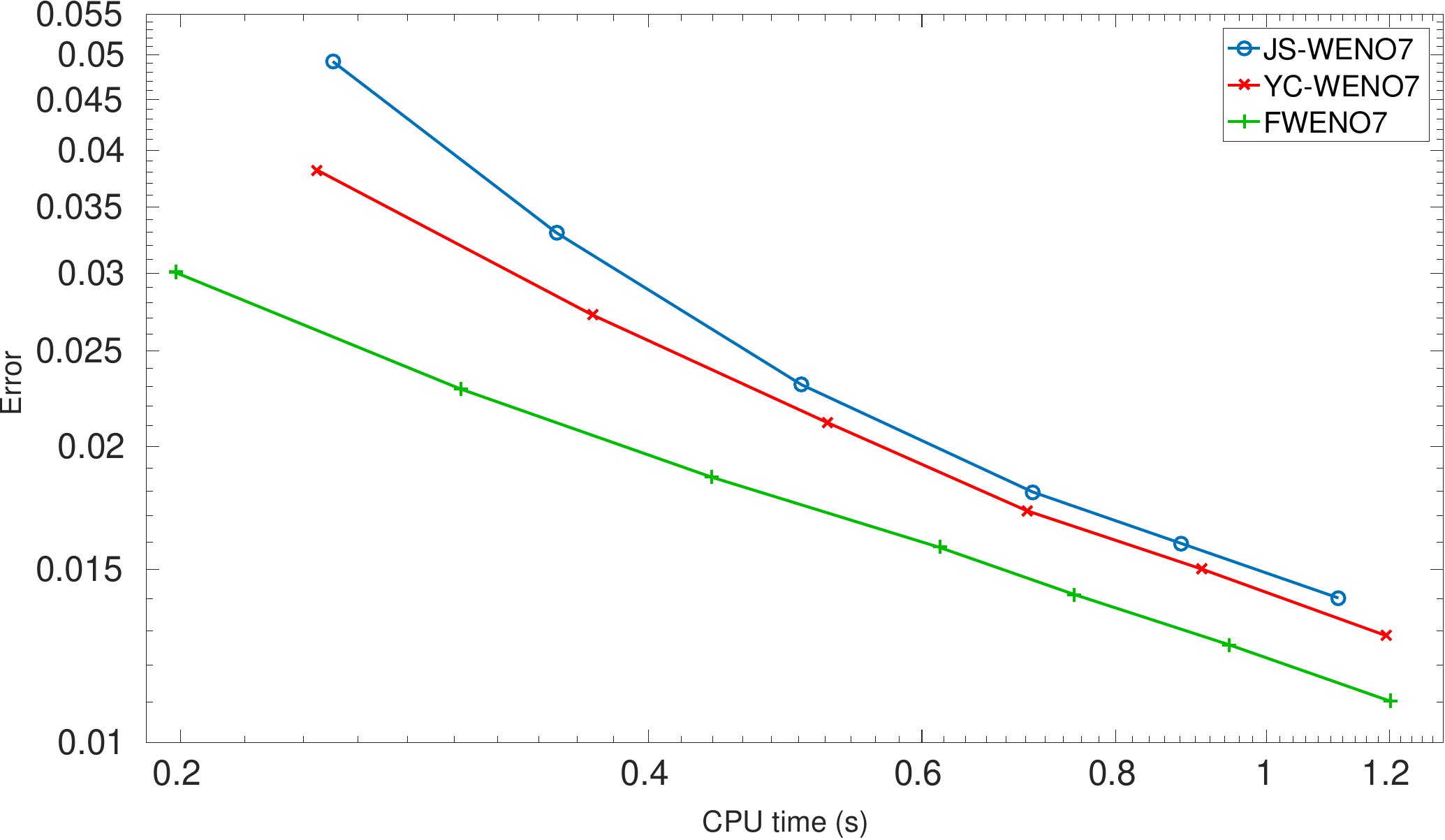}
  \end{tabular}
  \caption{Example 3: (Shu-Osher problem, 1D Euler equations of gas dynamics): efficiency plot for seventh-order schemes, corresponding to the 
   numerical solution at $T=1.8$.}
  \label{cpu_o7}
\end{figure}

\begin{figure}[t]
  \centering
  \begin{tabular}{cc}
    \includegraphics[width=0.85\textwidth]{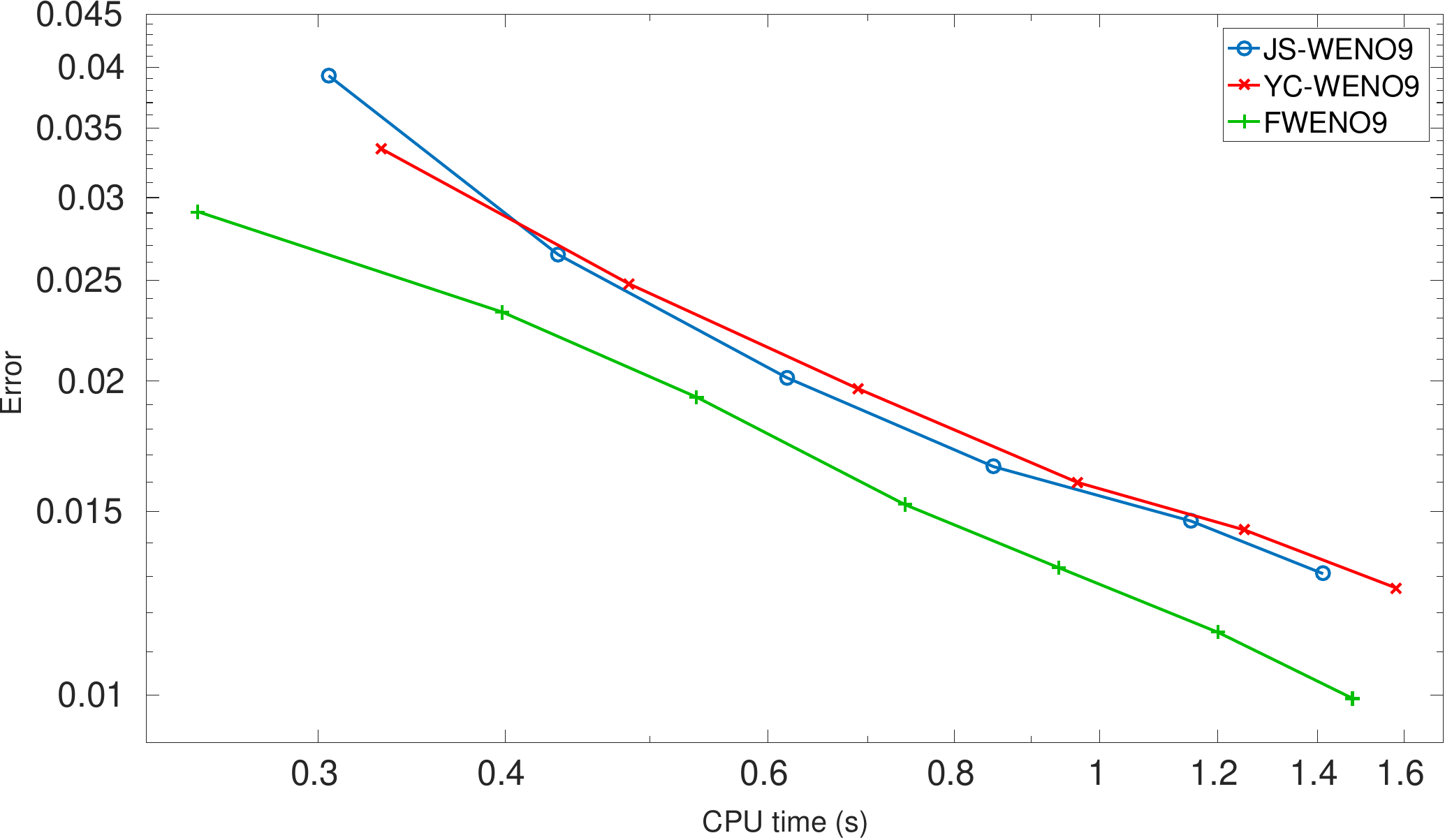}
  \end{tabular}
  \caption{Example 3: (Shu-Osher problem, 1D Euler equations of gas dynamics): efficiency plot for ninth-order schemes, corresponding to the 
   numerical solution at $T=1.8$.}
  \label{cpu_o9}
\end{figure}

The 1D Euler equations for gas dynamics are given by 
 $\boldsymbol{u} = ( \rho, \rho v, E)^{\mathrm{T}}$ and $\boldsymbol{f} (\boldsymbol{u}) = 
  \boldsymbol{f}^1  (\boldsymbol{u}) = (\rho v, 
        p+\rho v^2, 
        v(E+p))^{\mathrm{T}}$,  where $\rho$~is the density, $v$~is the velocity and  $E$~is the specific energy of the system. The variable~$p$ stands for the pressure and is given by the equation of state
$$p=\left(\gamma-1\right)\left(E-\frac{1}{2}\rho v^2\right),$$
where $\gamma$~is the adiabatic constant that will be taken as $\gamma =1.4$. The spatial domain is  $\Omega:=(-5,5)$, 
 and    the initial condition   
\begin{align*} 
(\rho,v,p) (x, 0) = 
\begin{cases} \displaystyle 
  \biggl(\frac{27}{7}, 
    \frac{4\sqrt{35}}{9}, 
    \frac{31}{3} \biggr) & \text{if  $x\leq-4$,}  \\[4mm]  \displaystyle
   \biggl(1+\frac{1}{5}\sin(5x), 0, 1 \biggr) 
    & \text{if $x>-4$,} 
\end{cases} \end{align*} 
stipulates 
 the interaction of a Mach~3 shock with a sine wave and is complemented  with left inflow and right outflow boundary conditions.
We run the simulation until $T=1.8$ and compare the schemes against a reference solution computed with a resolution of $N=16000$. 
Figures~\ref{shuosher_o5_n=200} and~\ref{shuosher_o5_n=400} display the results obtained by  fifth-order schemes with resolutions of $N=200$ and $N=400$ cells, respectively. 
  The results are similar for  the YC-WENO and FWENO schemes and are slightly less sharply resolved for the  JS-WENO scheme. In order to highlight the superior performance of the FWENO scheme, we also plot the numerical error against the CPU time for schemes of order $2r-1$, with $3\leq r\leq 5$, which is shown in Figures \ref{cpu_o5}--\ref{cpu_o9}. These results clearly indicate that 
 the scheme with the new smoothness indicators is more efficient (in terms of error reduction versus CPU time) than the schemes that employ the traditional Jiang-Shu smoothness indicators.

\subsubsection*{Example 4 (Sod shock tube problem,  1D Euler equations of gas dynamics)}

\begin{figure}[t]
  \centering
  \begin{tabular}{cc}
    \includegraphics[width=0.47\textwidth]{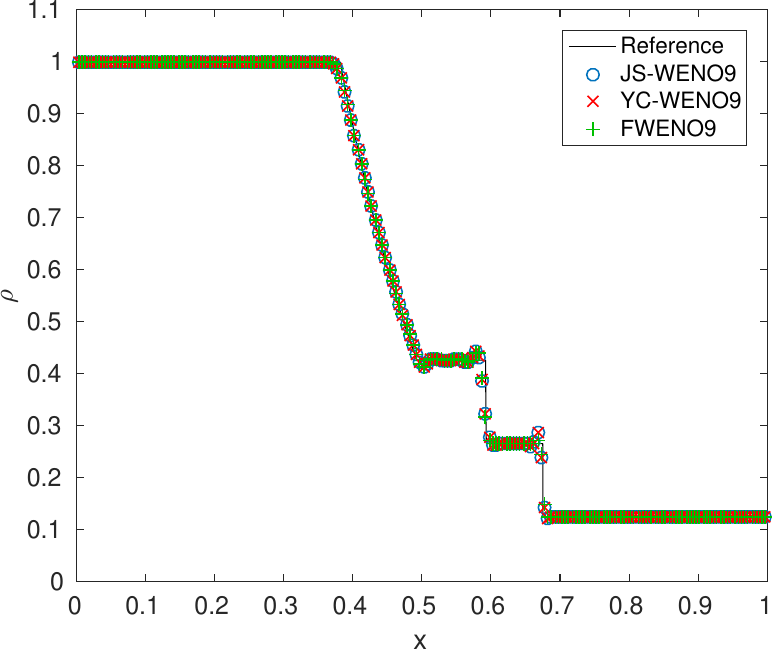} & \includegraphics[width=0.47\textwidth]{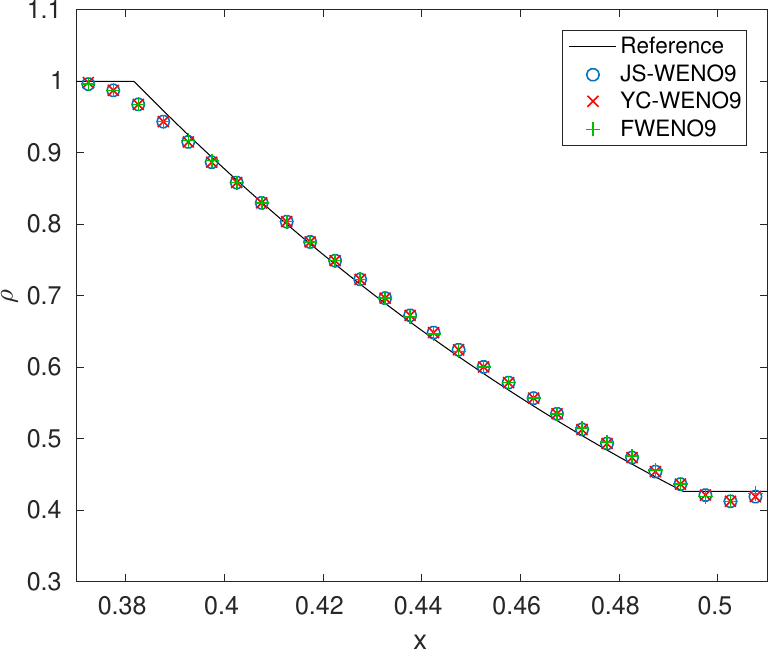} \\
    \includegraphics[width=0.47\textwidth]{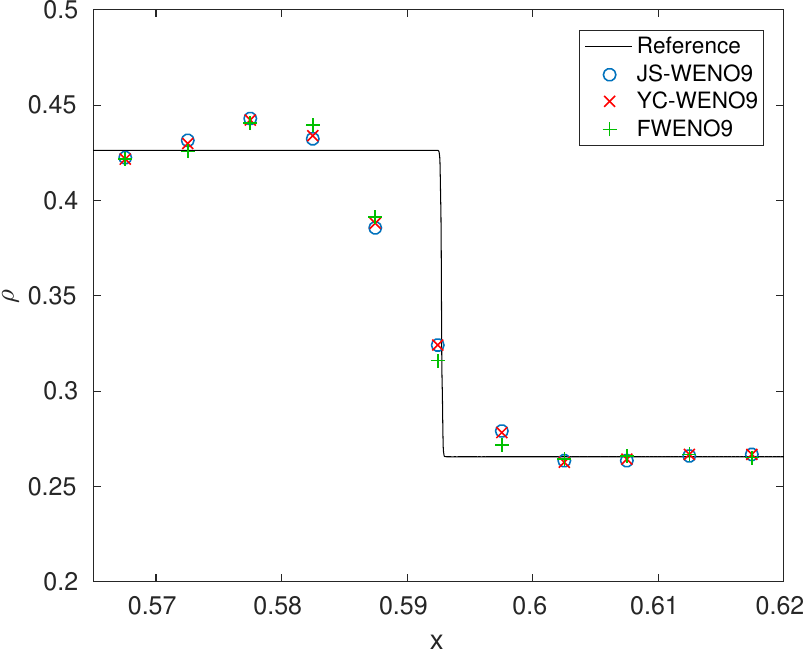} & \includegraphics[width=0.47\textwidth]{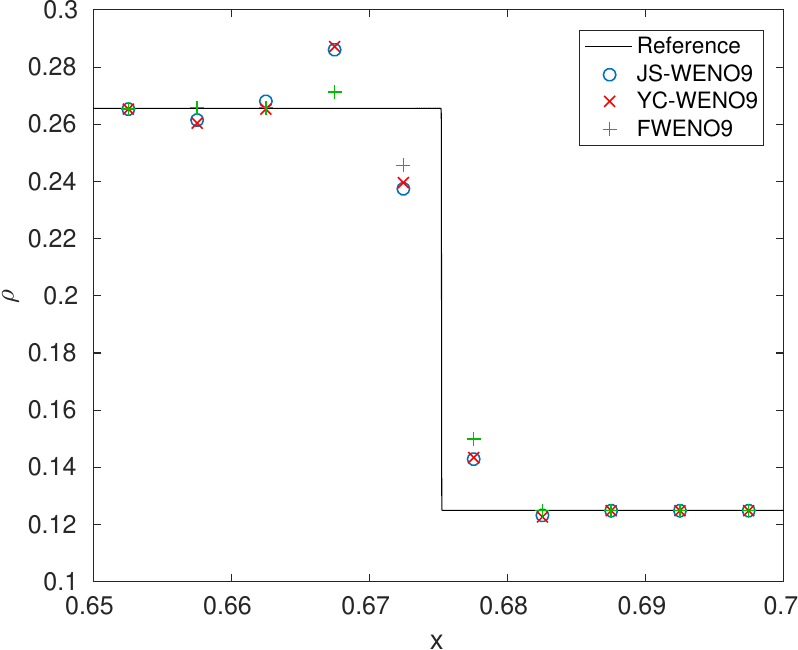}
  \end{tabular}
  \caption{Example 4 (Sod shock tube problem, 1D Euler equations of gas dynamics): numerical solution for $N=200$ (top left) and enlarged views (top right, bottom left and bottom right) at $T=0.1$ produced by 
   ninth-order schemes with conventional (JS-WENO9 and YC-WENO9) and 
   new (FWENO9) smoothness indicators. The reference solution with $N=100000$ cells is also shown.}
  \label{sod_o9_n=200}
\end{figure}

\begin{figure}[t] 
  \centering
  \begin{tabular}{cc}
    \includegraphics[width=0.85\textwidth]{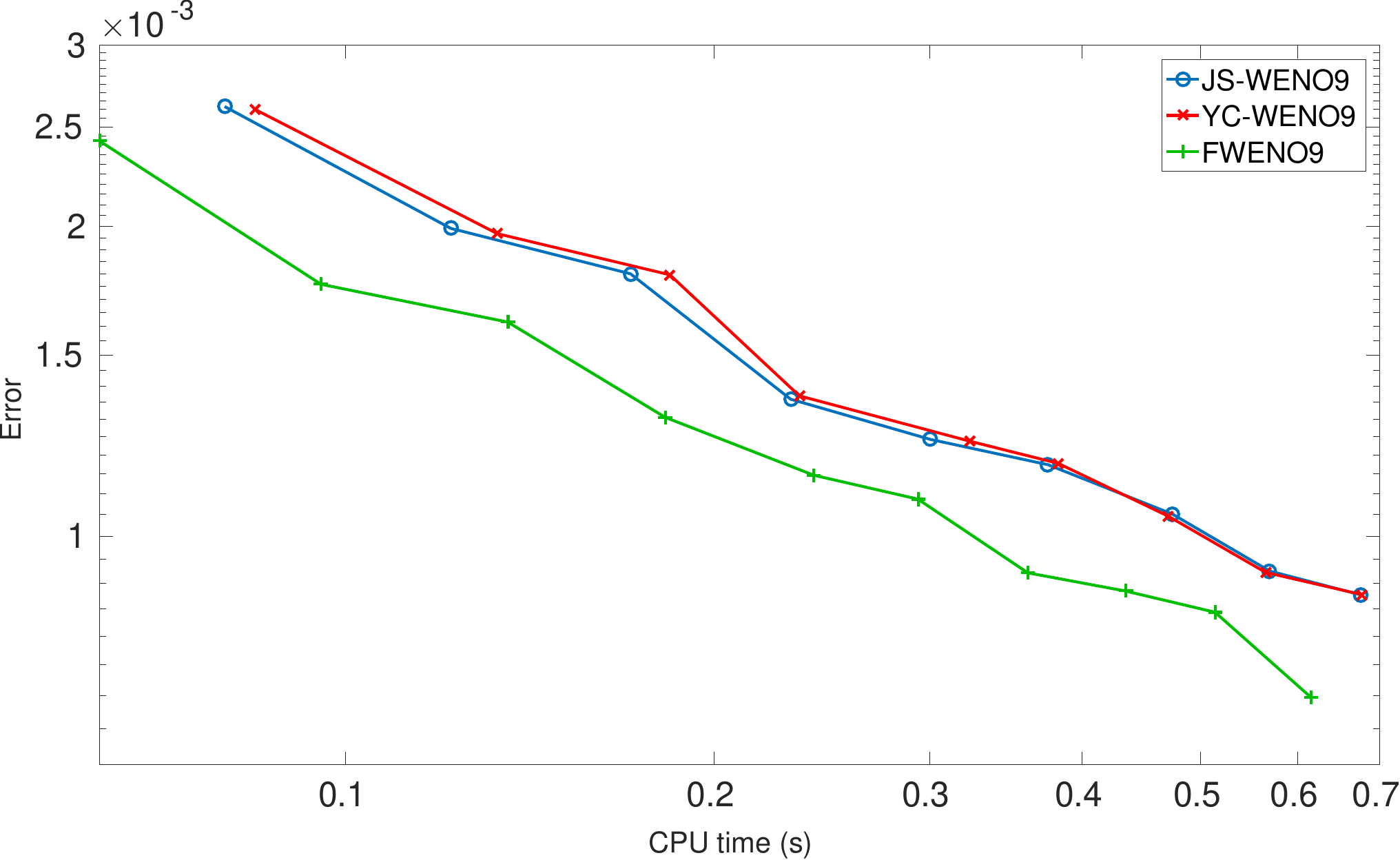}
  \end{tabular}
  \caption{Example 4 (Sod shock tube problem, 1D Euler equations of gas dynamics):  efficiency plot for ninth-order schemes, corresponding to the 
   numerical solution at $T=0.1$.}
  \label{effsod_o9}
\end{figure}

We now apply WENO schemes to the  1D Euler equations of gas dynamics on $\Omega=(0,1)$ with the initial condition
\begin{align*} 
(\rho,v,p) (x, 0) = 
\begin{cases}
  (1, 0, 1) & \text{if  $x\leq0.5$,}  \\
   (0.125, 0, 0.1) 
    & \text{if $x>0.5$} 
\end{cases} \end{align*}
and left and right Dirichlet boundary conditions corresponding to the shock tube problem proposed by Sod \cite{Sod1978}. This problem has been tackled in many other papers afterwards, such as in \cite{BurgerKumarZorio2017}. The numerical result is 
 produced by  ninth-order schemes with a resolution of $N=200$ cells compared against a reference solution which has been computed with a resolution of $N=100000$ by the classical JS-WENO scheme. The simulation is run until $T=0.1$ and the results are depicted in Figure~\ref{sod_o9_n=200}.
  It turns our that the schemes produce very  similar results. In fact, the most remarkable differences are favorable to our proposed FWENO scheme, 
   since the behavior is slightly  less oscillatory  near the contact discontinuity and the shock.
Finally, an efficiency comparison is presented in Figure \ref{effsod_o9}, where it can be concluded that our FWENO scheme turns out to be again more efficient than their classical JS-WENO and YC-WENO counterparts.

\subsection{2D conservation law experiments}\label{subsec:2d_cons_laws}

\subsubsection*{Example 5 (double Mach reflection, 2D Euler equations of gas dynamics)}

\begin{figure}[t]
  \centering
  \setlength{\tabcolsep}{1cm}
  \begin{tabular}{cc}
    \includegraphics[width=0.38\textwidth]{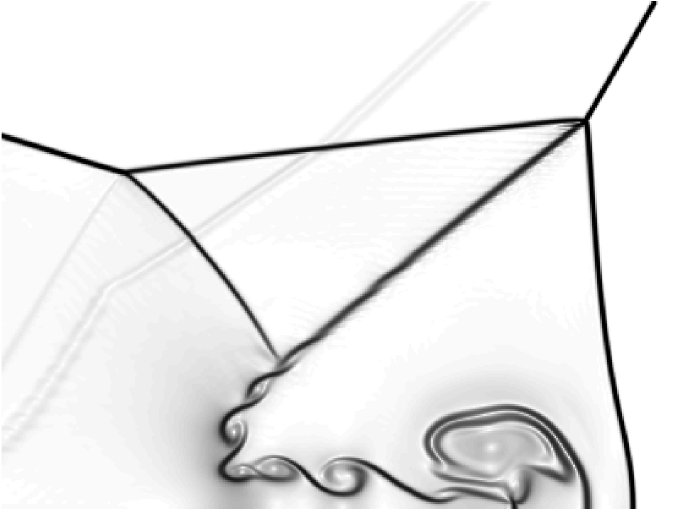} & \includegraphics[width=0.38\textwidth]{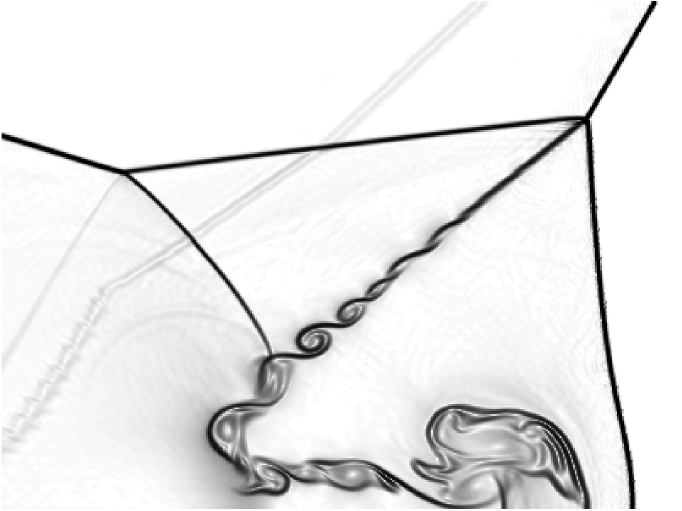} \\
    (a) JS-WENO5, $2048\times512$ & (b) YC-WENO5, $2048\times512$ \\
    \includegraphics[width=0.38\textwidth]{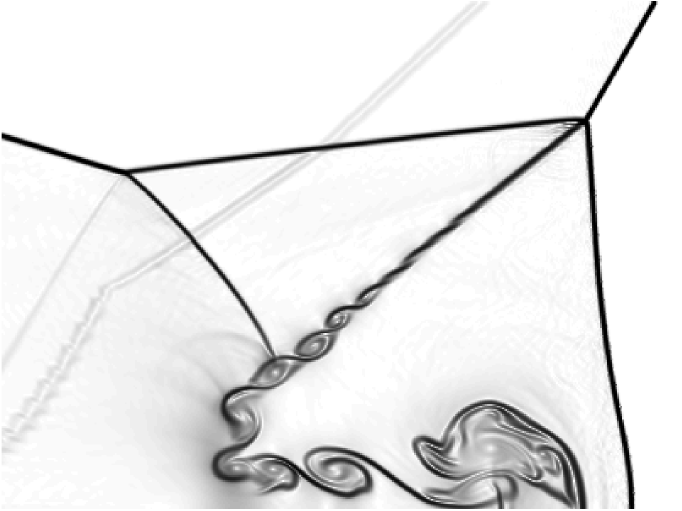} & \includegraphics[width=0.38\textwidth]{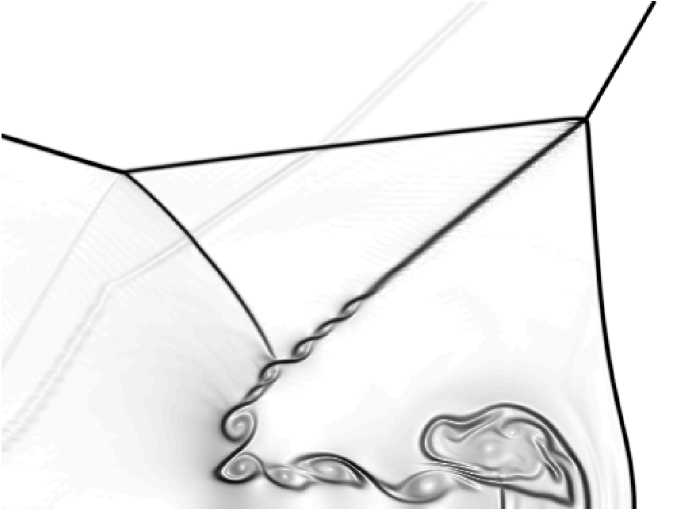} \\
    (c) FWENO5, $2048\times512$ & (d) JS-WENO5, $2560\times640$ \\
    \includegraphics[width=0.38\textwidth]{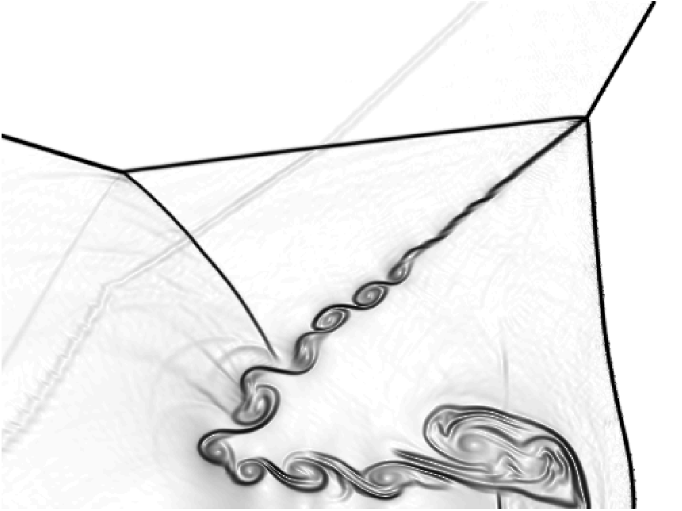} & \includegraphics[width=0.38\textwidth]{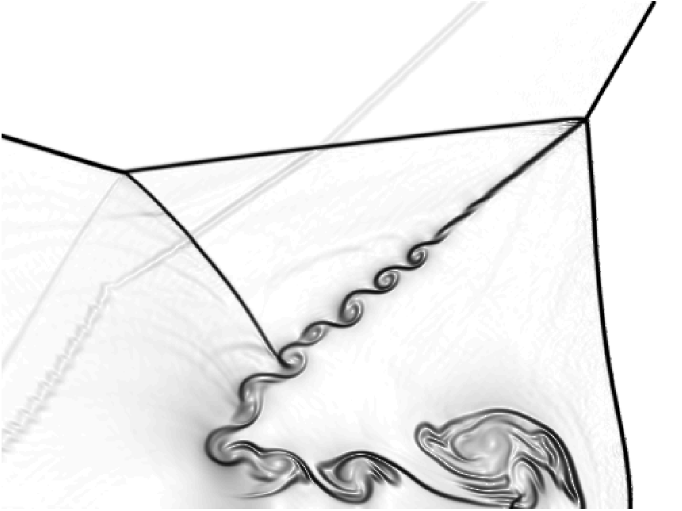}  \\
    (e) YC-WENO5, $2560\times640$ & (f) FWENO5, $2560\times640$ \\
  \end{tabular}
  \caption{Example~5 (double Mach reflection, 2D Euler equations of gas dynamics): enlarged views  of turbulent zone of numerical solutions at $T=0.2$ (Schlieren plot).}
  \label{dmr}
\end{figure}

The two-dimensional Euler equations for inviscid gas dynamics
are given by
$$\boldsymbol{u}_t+\boldsymbol{f}^1(\boldsymbol{u})_x+\boldsymbol{f}^2(\boldsymbol{u})_y=0,$$
 with  \begin{align*} 
& \boldsymbol{u} =  \begin{pmatrix} 
    \rho \\
    \rho v^x \\
    \rho v^y \\
    E \end{pmatrix}, \quad 
  \boldsymbol{f}^1 (\boldsymbol{u})= \begin{pmatrix}  
    \rho v^x \\
    p+\rho (v^x)^2 \\
    \rho v^xv^y \\
    v^x(E+p) 
  \end{pmatrix},  \quad   \boldsymbol{f}^2 (\boldsymbol{u})= \begin{pmatrix} 
    \rho v^y \\
    \rho v^xv^y \\
    p+\rho (v^y)^2 \\
    v^y(E+p) 
  \end{pmatrix}. 
\end{align*}
Here  $\rho$ is the density, $(v^x, v^y)$  is the velocity, $E$ is the specific energy, and  $p$  is the pressure  that  is given by the equation of state $$p=(\gamma-1)\left(E-\frac{1}{2}\rho((v^x)^2+(v^y)^2)\right),$$ where
 the adiabatic constant is again chosen as  $\gamma =1.4$.
This experiment uses these equations to model a vertical right-going Mach 10 shock colliding with an equilateral triangle. By symmetry, this is equivalent to a collision with a ramp with a slope of $30^{\circ}$ with respect to the horizontal line.

For  sake of simplicity, we consider the equivalent problem in a rectangle, consisting in a rotated shock, whose vertical angle is 
$30^{\circ}$. The domain is the rectangle $\Omega=[0,4]\times[0,1]$, whose initial conditions are
\begin{gather*}  (\rho,v^x,v^y,E) (x,y,0)=\begin{cases}
  \boldsymbol{c}_1= (\rho_1,v_1^x,v_1^y,E_1)    & \text{if $y\leq 1/4 +\tan(\pi/6)x$,} \\
  \boldsymbol{c}_2=  (\rho_2,v_2^x,v_2^y,E_2)     & \text{if $y > 1/4 +\tan(\pi/6)x$,}  
\end{cases} \\
    \boldsymbol{c}_1 =
    \bigl(8,8.25\cos(\pi/6),-8.25\sin(\pi/6),563.5\bigr), \quad 
    \boldsymbol{c}_2=    (1.4,0,0,2.5).
  \end{gather*} 
We impose inflow boundary conditions, with value $\boldsymbol{c}_1$, at the left side, $\{0\}\times[0,1]$, outflow boundary conditions both at $[0,1/4]\times\{0\}$ and $\{4\}\times[0,1]$, reflecting boundary conditions at  $(1/4,4]\times\{0\}$ and inflow boundary conditions at the upper side, $[0,4]\times\{1\}$, which mimics the shock at its actual traveling speed:
\begin{align*} 
 (\rho,v^x,v^y,E) (x,1,t)=\begin{cases}
  \boldsymbol{c}_1 & \text{if $x\leq 1/4 + (1+20t)/\sqrt{3}$,}  \\
  \boldsymbol{c}_2 & \text{if $x>1/4 + (1+20t)/\sqrt{3}$.}  
\end{cases} \end{align*}
We run different simulations until $T=0.2$ both at a resolution of $2048\times512$ points and a resolution of $2560\times640$ points, shown in 
Figure~\ref{dmr}, in both cases with $\textnormal{CFL}=0.4$ and involving the JS-WENO scheme, the YC-WENO method and our FWENO scheme for the case of fifth-order accuracy. 
The results show that both YC-WENO and FWENO schemes produce sharper resolution than JS-WENO, and in turn they have similar resolution between then.

\subsubsection*{Example 6: Riemann problem}

\begin{figure}[t]
  \centering
  \setlength{\tabcolsep}{1cm}
  \begin{tabular}{cc}
    \includegraphics[width=0.38\textwidth]{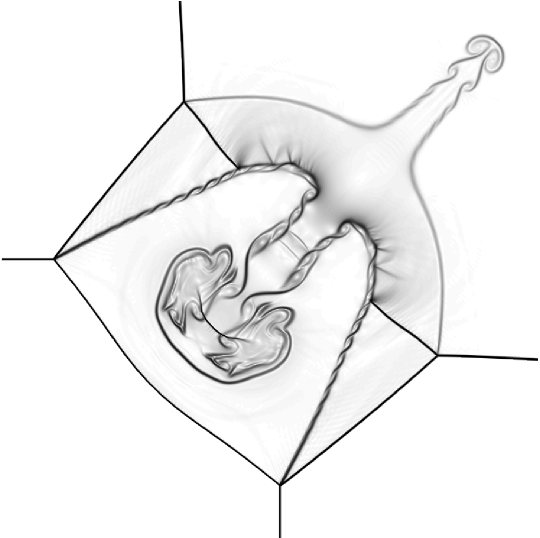} & \includegraphics[width=0.38\textwidth]{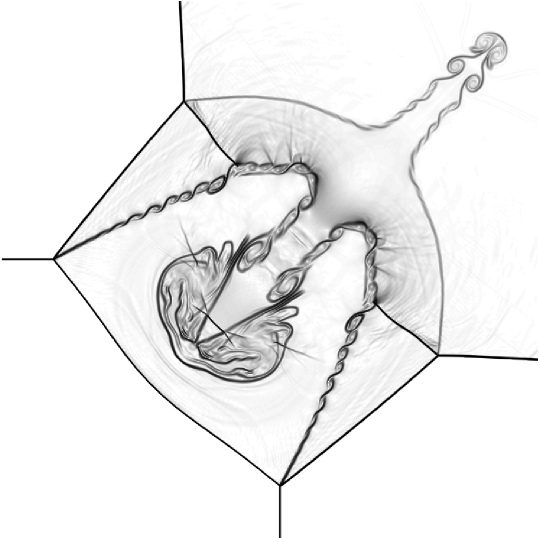} \\
    (a) JS-WENO5, $2048\times2048$ & (b) YC-WENO5, $2048\times2048$ \\
    \includegraphics[width=0.38\textwidth]{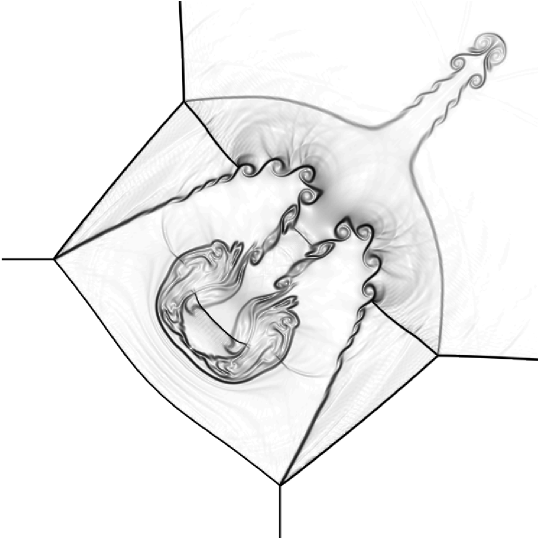} & \includegraphics[width=0.38\textwidth]{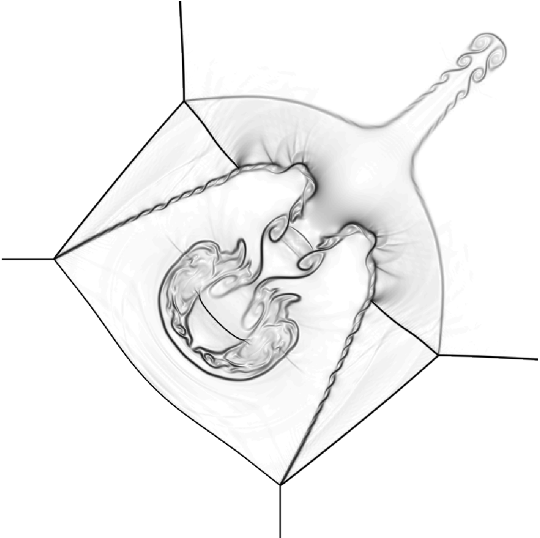} \\
    (c) FWENO5, $2048\times2048$ & (d) JS-WENO5, $2560\times2560$ \\
    \includegraphics[width=0.38\textwidth]{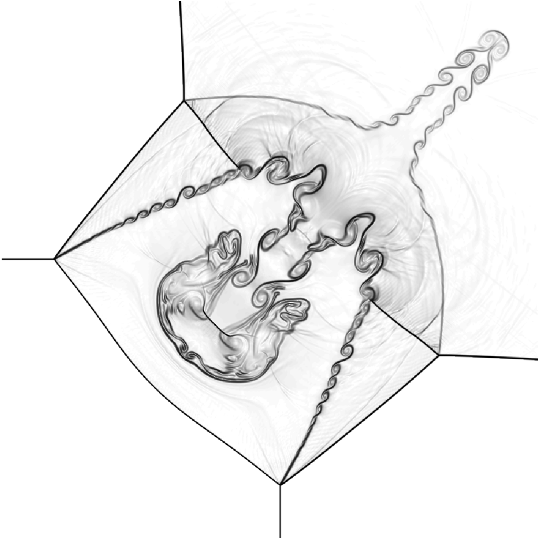} & \includegraphics[width=0.38\textwidth]{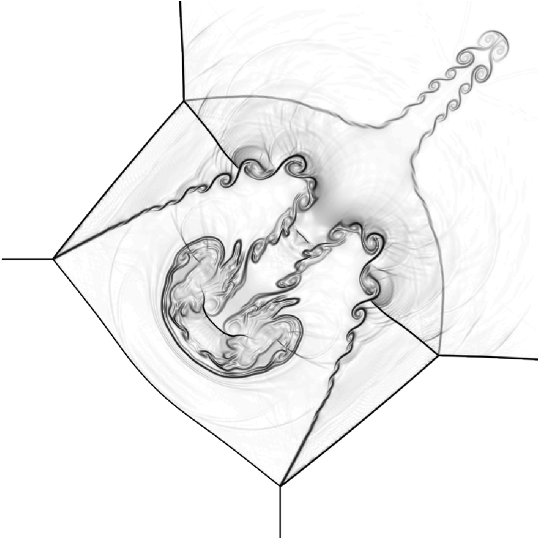} \\
    (e) YC-WENO5, $2560\times2560$ & (f) FWENO5, $2560\times2560$
  \end{tabular}
  \caption{Example~6 (Riemann problem, 2D Euler equations of gas dynamics): enlarged views  of turbulent zone of numerical solutions at $T=0.3$ (Schlieren plot).}
  \label{riemann}
\end{figure}

Finally, we solve numerically a Riemann problem for the 2D Euler equations on the domain $(0,1)\times(0,1)$. Riemann problems for 2D Euler equations were first studied in \cite{SchulzRinne}. The initial data is taken from \cite[Sect.~3, Config.~3]{KurganovTadmor}:
\begin{align*}
  \boldsymbol{u}(x, y, 0)=(\rho(x, y, 0), \rho(x, y, 0)v^x(x, y, 0), \rho(x, y, 0)v^y(x, y, 0), E(x, y, 0))
  \end{align*}
and
\[
\begin{pmatrix}
  \rho(x, y, 0)\\
  v^x(x, y, 0)\\
  v^y(x, y, 0)\\
  p(x, y, 0)
\end{pmatrix}^{\mathrm{T}} 
  =
\begin{cases}
  (1.5, 0, 0, 1.5) & \text{for $x>0.5$, $y>0.5$,}  \\
  (0.5323, 1.206, 0, 0.3) & \text{for  $x\leq0.5$, $y>0.5$,}  \\
  (0.138, 1.206, 1.206, 0.029) &  \text{for $x\leq0.5$, $y\leq0.5$,}  \\
  (0.5323, 0, 1.206, 0.3) &  \text{for $x>0.5$, $y\leq0.5$,} 
\end{cases}
\]
with the same equation of state as in the previous test.
The simulation is performed taking $s_2=2$, with the final time $T=0.3$, $\textnormal{CFL}=0.4$, resolutions $2048\times2048$ and $2560\times2560$ and comparing the same schemes with the same parameters as in Example~5. The results are shown in Figure \ref{riemann}. 
The same conclusions as in Example~5 are drawn.

\begin{figure}[t] 
  \centering
  \begin{tabular}{c}
    \includegraphics[width=0.715\textwidth]{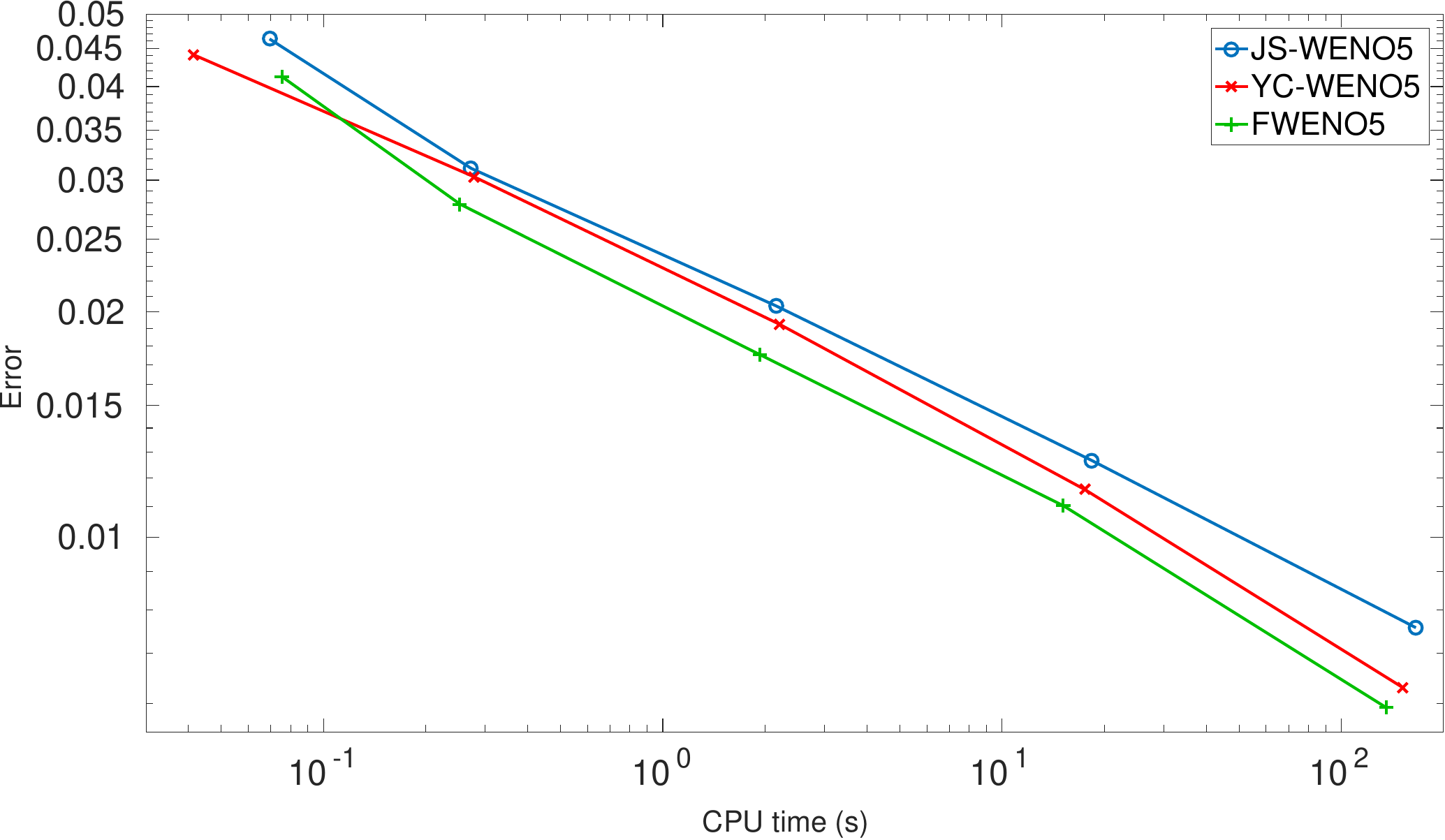}\\
    \includegraphics[width=0.715\textwidth]{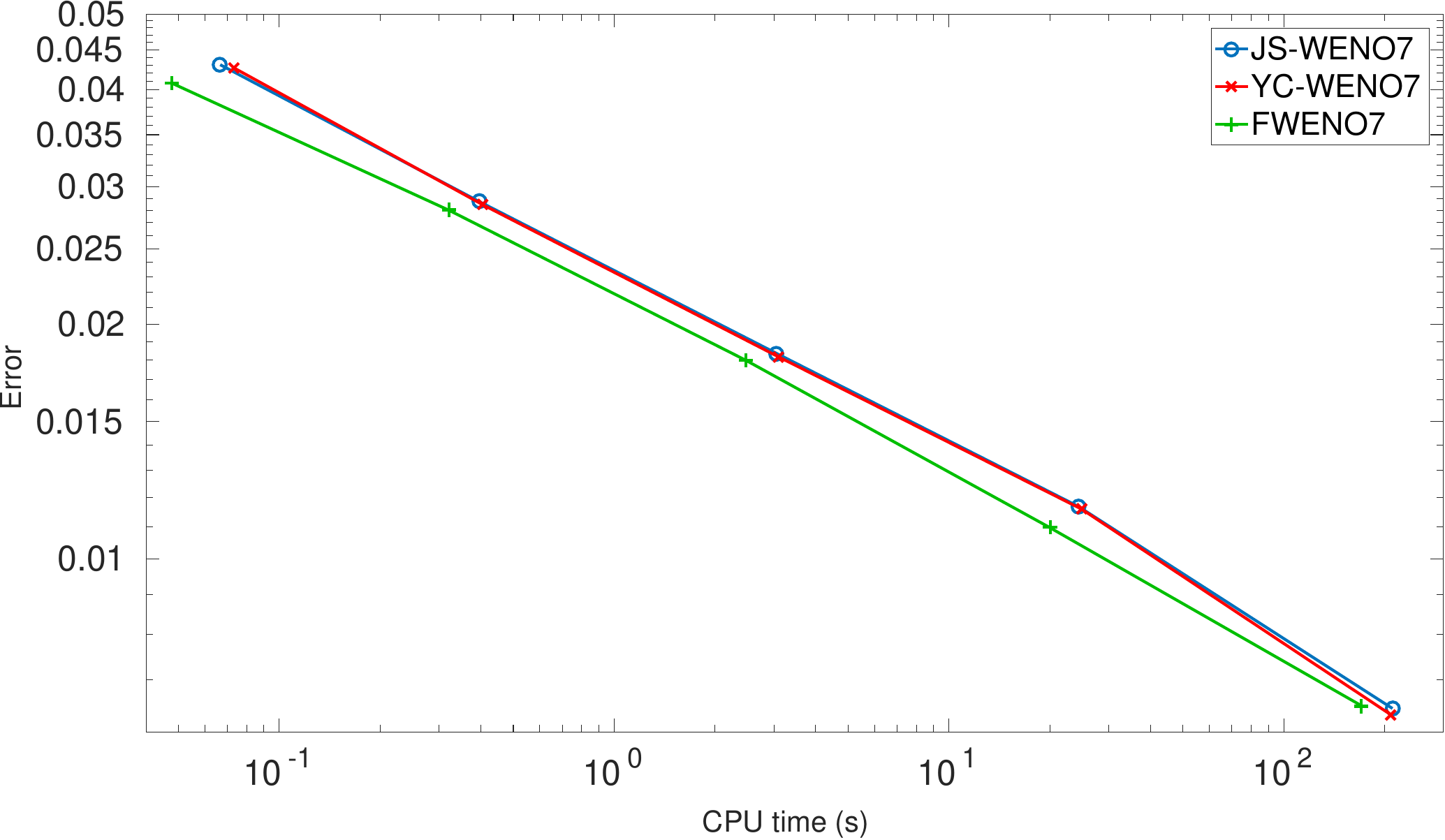}\\
    \includegraphics[width=0.715\textwidth]{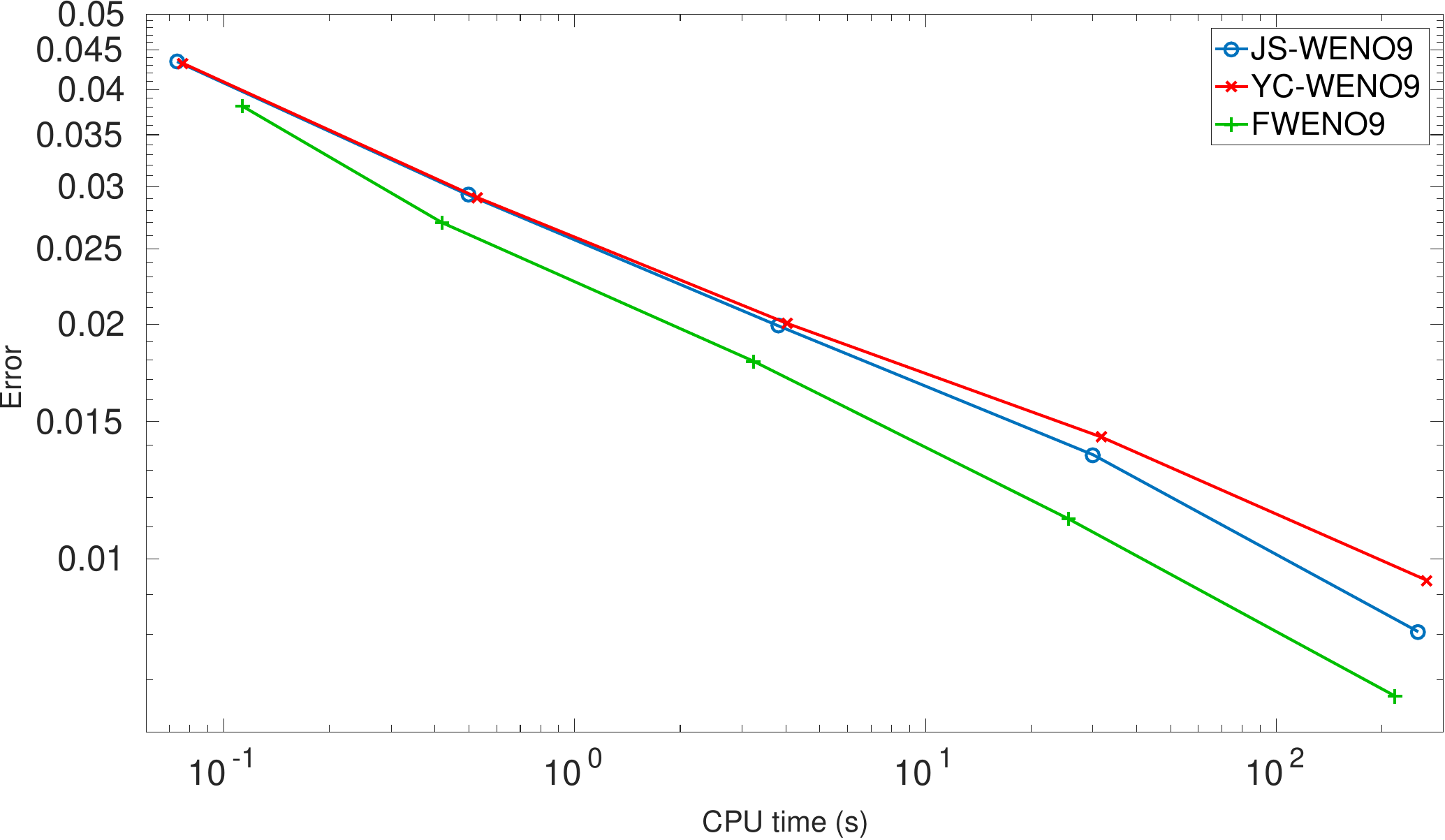}
  \end{tabular}
  \caption{Example 6 (2D Riemann problem, 2D Euler equations of gas dynamics):  efficiency plot for fifth, seventh and ninth-order schemes, corresponding to the 
   numerical solution at $T=0.3$.}
  \label{effriemann}
\end{figure}

We now use the solutions computed with the grid of $2560\times2560$ points as reference solutions to perform efficiency tests by comparing error versus CPU time involving numerical solutions with grid sizes $16\cdot2^n\times16\cdot2^n$, $0\leq n\leq 4$, for the corresponding fifth, seventh and ninth order schemes. The results are shown in Figure \ref{effriemann} and again indicate a higher performance for the FWENO scheme.

\section{Conclusions}
\label{sec:conclusions}

In this paper a set of alternative smoothness indicators, cheaper than the classical Jiang-Shu ones, has been presented. The theoretical results show that when used in Yamaleev-Carpenter type weight constructions they attain the same accuracy properties than the ones obtained with the original smoothness indicators. The numerical experiments confirm all these theoretical considerations. Also, the numerical evidence obtained in the problems from hyperbolic conservation laws with weak solutions also shows that the quality of the approximation is similar in both cases, being the schemes with the modified smoothness indicators (FWENO) more efficient than their traditional counterparts (JS-WENO and YC-WENO).

As for future work, we encompass extrapolating the benefits of this new weight design with the simplified smoothness indicators in the context of WENO extrapolation for numerical boundary conditions and even generalized WENO interpolations/extrapolations in the context of non-uniform grids, in which the computational benefits of using these new smoothness indicators with respect to the traditional ones are expected to be much higher than in uniform grids.

\section*{Acknowledgements} 

AB, PM and DZ are supported by Spanish MINECO project
MTM2017-83942-P.  
RB is supported by  CONICYT/PIA/Concurso Apoyo a Centros Cient\'{\i}ficos y Tecnol\'{o}gicos de Excelencia con Financiamiento Basal AFB170001; Fondecyt project 1170473; and CRHIAM, project CONICYT/FONDAP/15130015. 
PM is also
supported by Conicyt (Chile),   project PAI-MEC, folio 80150006.  
 DZ is also supported by Conicyt (Chile) through Fondecyt project 3170077.


\begin{thebibliography}{10}
  \providecommand{\url}[1]{{#1}}
  \providecommand{\urlprefix}{URL }
  \expandafter\ifx\csname urlstyle\endcsname\relax
  \providecommand{\doi}[1]{DOI~\discretionary{}{}{}#1}\else
  \providecommand{\doi}{DOI~\discretionary{}{}{}\begingroup
    \urlstyle{rm}\Url}\fi

\bibitem{SINUM2011}
  Aràndiga, F., Baeza, A., Belda, A.M., Mulet, P.: Analysis of {WENO} schemes
  for full and global accuracy.
  \newblock SIAM J.\ Numer.\ Anal. \textbf{49}(2), 893--915 (2011)

\bibitem{BurgerKumarZorio2017}
  Bürger, R., Kenettinkara, S. K., Zorío, D.: Approximate Lax-Wendroff discontinuous Galerkin methods for hyperbolic conservation laws.
  \newblock Computers \& Mathematics with Applications, \textbf{74}(6), 1288--1310 (2017)

\bibitem{DonatMarquina96}
  Donat, R., Marquina, A.: Capturing shock reflections: An improved flux formula.
  \newblock J. Comput. Phys. \textbf{125}, 42--58 (1996)

\bibitem{JiangShu96}
  Jiang, G.S., Shu, C.W.: Efficient implementation of {Weighted} {ENO} schemes.
  \newblock J. Of Comput. Phys. \textbf{126}, 202--228 (1996)

\bibitem{KurganovTadmor}
  Kurganov, A., Tadmor, E.: Solution of two-dimensional {R}iemann problems for gas dynamics without Riemann problem solvers
  \newblock Numer.\ Methods  Partial Differential Equations  \textbf{18}, 584--608 (2002)

\bibitem{LiuOsherChan94}
  Liu, X-D., Osher, S., Chan, T.: Weighted essentially non-oscillatory schemes.
  \newblock J. Comput. Phys. \textbf{115}(1) 200--212 (1994)

\bibitem{SchulzRinne}
  Schulz-Rinne, C. W.: Classification of the Riemann problem for two-dimensional gas dynamics.
  \newblock SIAM Journal on Mathematical Analysis \textbf{24}(1), 76--88 (1993)

\bibitem{ShuOsher89}
  Shu, C.-W., Osher, S.: Efficient implementation of essentially non-oscillatory
  shock-capturing schemes.
  \newblock J.\ Comput.\ Phys. \textbf{77}, 439--471 (1988)

\bibitem{ShuOsher1989}
  Shu, C.-W., Osher, S.: Efficient implementation of essentially non-oscillatory
  shock-capturing schemes, {II}.
  \newblock J. Comput. Phys. \textbf{83}(1), 32--78 (1989)

\bibitem{shu98} 
  Shu, C.-W.: Essentially non-oscillatory and weighted essentially non-oscillatory schemes for hyperbolic conservation laws. In 
  Cockburn, B., Johnson, C., Shu, C.-W. and Tadmor, E.  (Quarteroni, A., Ed.):  Advanced Numerical Approximation of Nonlinear Hyperbolic Equations. 
  \newblock Lecture Notes in Mathematics, vol.~1697, Springer-Verlag, Berlin, 325--432 (1998) 

\bibitem{Sod1978}
  Sod, G. A.: A survey of several finite difference methods for systems of nonlinear hyperbolic conservation laws.
  \newblock J. Comput. Phys. \textbf{27}(1), 1--31 (1978)

\bibitem{YamaleevCarpenter2009O3}
  Yamaleev, N. K., Carpenter, M. H.: Third-order Energy Stable WENO scheme.
  \newblock J. Comput. Phys. \textbf{228}(8), 3025--3047 (2009)

\bibitem{YamaleevCarpenter2009}
  Yamaleev, N. K., Carpenter, M. H.: A systematic methodology to for constructing high-order energy stable WENO schemes.
  \newblock J. Comput. Phys. \textbf{228}(11), 4248--4272 (2009)
  
\end{thebibliography}
\end{document}